\documentclass[12pt]{article}
\pdfoutput=1

\usepackage{amsmath,amssymb,amsthm,a4wide}
\usepackage{tikz}
\usepackage{mathtools}

\newcommand{\Z}{\mathbb{Z}}

\newcommand{\C}{\mathbb{C}}

\newcommand{\f}{f}
\newcommand{\Q}[1]{{Q}^{({#1})}}
\newcommand{\q}[1]{q^{({#1})}}

\newcommand{\osp}{\mathfrak{osp}}

\newcommand{\gl}{\mathfrak{gl}}

\newcommand{\B}{\mathcal{B}}
\newcommand{\BB}{\hat{\mathcal{B}}}

\newcommand{\e}[1]{\left\langle {#1} \right\rangle}

\newcommand{\str}{\mathop{\rm str}\nolimits}

\newcommand{\End}{\mathop{\rm End}\nolimits}

\newtheorem{theorem}{Theorem}[section]
\newtheorem{proposition}[theorem]{Proposition}
\newtheorem{lemma}[theorem]{Lemma}
\newtheorem{corollary}[theorem]{Corollary}
\newtheorem{conjecture}[theorem]{Conjecture}
\newtheorem{definition}[theorem]{Definition}

\theoremstyle{remark}
\newtheorem{remark}[theorem]{Remark}

\newcommand{\aff}{{\rm aff}}
\newcommand{\cri}{{\rm cri}}

\newcommand{\non}{\nonumber}

\newcommand{\wh}{\widehat}

\newcommand{\ot}{\otimes}
\newcommand{\la}{\lambda}

\newcommand{\ga}{\gamma}

\newcommand{\ep}{\epsilon}
\newcommand{\eee}{\epsilon^{}}
\newcommand{\ka}{\kappa}

\newcommand{\ta}{\ve}

\newcommand{\vp}{\varphi}

\newcommand{\de}{\delta}

\newcommand{\om}{\omega}

\newcommand{\su}{s^{}}
\newcommand{\hra}{\hookrightarrow}
\newcommand{\ve}{\varepsilon}
\newcommand{\ts}{\,}

\newcommand{\tss}{\hspace{1pt}}
\newcommand{\pr}{^{\tss\prime}}

\newcommand{\U}{ {\rm U}}

\newcommand{\CC}{\mathbb{C}\tss}

\newcommand{\ZZ}{\mathbb{Z}\tss}

\newcommand{\Dc}{\mathcal{D}}

\newcommand{\Bc}{\mathcal{B}}

\newcommand{\Yc}{\mathcal{Y}}

\newcommand{\ev}{{\rm ev}}

\newcommand{\oa}{\mathfrak{o}}
\newcommand{\spa}{\mathfrak{sp}}
\newcommand{\g}{\mathfrak{g}}

\newcommand{\z}{\mathfrak{z}}

\newcommand{\sll}{\mathfrak{sl}}

\newcommand{\bi}{\bar{\imath}}
\newcommand{\bj}{\bar{\jmath}}
\newcommand{\bk}{\bar{k}}
\newcommand{\bl}{\bar{l}}

\newcommand{\Sym}{\mathfrak S}

\newcommand{\Fand}{\qquad\text{and}\qquad}

\numberwithin{equation}{section}

\newcommand{\bal}{\begin{aligned}}
\newcommand{\eal}{\end{aligned}}
\newcommand{\beq}{\begin{equation}}
\newcommand{\eeq}{\end{equation}}
\newcommand{\ben}{\begin{equation*}}
\newcommand{\een}{\end{equation*}}

\newcommand{\bpf}{\begin{proof}}
\newcommand{\epf}{\end{proof}}

\def\beql#1{\begin{equation}\label{#1}}

\usepackage{hyperref}
\hypersetup{colorlinks,citecolor=blue,filecolor=black,linkcolor=blue,urlcolor=blue}
\usepackage{enumerate}

\usepackage{lmodern}     
\usepackage[T1]{fontenc}

\begin{document}

\let\savedleq=\leq 
\let\savedgeq=\geq %
\let\leq=\leqslant 
\let\geq=\geqslant %

\title{\Large\bf Segal--Sugawara vectors for orthosymplectic\\ Lie superalgebras}

\author{{Alexander Molev\quad and\quad Madeline Nurcombe}}

\date{} 
\maketitle


\begin{abstract}
We consider the centre of
the affine vertex algebra at the critical level
associated with the orthosymplectic Lie superalgebra.
It is well-known that the centre is a commutative superalgebra, and we construct
a family of its elements in an explicit form. In particular, this gives
a new proof of the formulas for the central elements for the orthogonal and symplectic
Lie algebras. Our arguments rely on the properties of a new extended Brauer-type algebra.

\end{abstract}

\section{Introduction}
\label{sec:int}

Let $\g$ be a finite-dimensional
Lie superalgebra over $\CC$ that is equipped with
an invariant super-symmetric bilinear form. Consider the corresponding
affine Kac--Moody superalgebra $\wh\g$ defined as a central extension
\ben
\wh\g=\g\tss[t,t^{-1}]\oplus\CC K
\een
of the Lie
superalgebra of Laurent polynomials $\g\tss[t,t^{-1}]$.
The vacuum module $V_{\ka}(\g)$ at the level $\ka\in\CC$
is the quotient of the
universal enveloping algebra $\U(\wh\g)$ by the left
ideal generated by $\g[t]$ and $K-\ka$.
The vacuum module has a vertex algebra structure and is known as
the {\em {\em({\em universal})} affine vertex algebra}; see e.g.
\cite{f:lc} and \cite{k:va} for definitions. The {\em centre}
of this vertex algebra is a commutative associative superalgebra
which can be regarded as a subalgebra of $\U(t^{-1}\g\tss[t^{-1}])$.

In the case of a simple Lie algebra $\g$, the centre is trivial, except at
the critical level $\ka=-h^{\vee}$, where $h^{\vee}$ is the dual Coxeter number for $\g$.
The vertex algebra $V_{-h^{\vee}}(\g)$ has a large centre $\z(\wh\g)$
which can be described by
\beql{ffcent}
\z(\wh\g)
=\{S\in V_{-h^{\vee}}(\g)\ |\ \g[t]\ts S=0\}.
\eeq
Any element of $\z(\wh\g)$
is called a {\em Segal--Sugawara vector\/}.
The algebra $\z(\wh\g)$ is equipped with the derivation $\tau=-d/dt$
arising from the vertex algebra structure.
By the celebrated theorem of Feigin and Frenkel~\cite{ff:ak}, the differential algebra
$\z(\wh\g)$ possesses generators
$S_1,\dots,S_n$ so that
$\z(\wh\g)$ is the algebra of polynomials
\ben
\z(\wh\g)=\CC[\tau^{\tss r}S_l\ |\ l=1,\dots,n,\ \ r\geqslant 0],
\een
where $n=\text{rank}\ts\g$; see also \cite{f:lc}. The algebra $\z(\wh\g)$
is known as the {\em Feigin--Frenkel centre}, and the elements $S_1,\dots,S_n$
form a {\em complete set of Segal--Sugawara vectors}.

Explicit formulas for complete sets of Segal--Sugawara vectors
were given in \cite{cm:ho} and \cite{ct:qs} for the Lie algebras $\g$ of type $A$,
and in \cite{m:ff} for types $B$, $C$ and $D$ with the use of the Brauer algebra.
A detailed exposition of these results, together with applications to
commutative subalgebras in enveloping algebras and to higher order Hamiltonians
in the Gaudin models, can be found in \cite{m:so}.
A complete set of Segal--Sugawara vectors
for the Lie algebra of type $G_2$ was produced in \cite{mrr:ss}
using computer-assisted calculations.
A different method to construct generators of $\z(\wh\g)$ was developed
in \cite{o:sf} which led to new explicit formulas in the case of the Lie algebras
of types $B,C,D$ and $G_2$. It was shown in \cite{m:ss} that they
coincide with those in \cite{m:ff} in the classical types.
Note also the recent work \cite{r:iff} where interpolating families of Segal--Sugawara operators
were constructed in the context of categorical vertex algebra theory.

Although no general analogue of the Feigin--Frenkel theorem is known
in the super case,
a few families of Segal--Sugawara vectors for the general linear Lie superalgebra
$\gl_{m|n}$ were constructed in \cite{mr:mm}, and it was conjectured therein that
these vectors generate the centre $\z(\wh\gl_{m|n})$ of the associated
affine vertex algebra at the critical level.
The conjecture
was proved for $\gl_{1|1}$ in \cite{mm:iv}, and for $\gl_{2|1}$ in \cite{an:ca}.
The affine vertex algebra associated with the orthosymplectic
Lie superalgebra $\osp_{1|2}$
was investigated in \cite{a:rs}. A conjectural description of its centre at the critical level
was pointed out in \cite[Remark~10]{a:rs}.

Our goal in this paper is to construct a family of Segal--Sugawara vectors for the
Lie superalgebra $\osp_{M|2n}$. As with the case of the orthogonal
and symplectic Lie algebras considered in \cite{m:ff}
(see also \cite[Ch.~8]{m:so}), we rely on the properties of a distinguished element $s^{(m)}$
(the {\em symmetriser})
of the Brauer algebra $\B_m(\om)$. According to the super version of the Schur--Weyl duality,
the actions of $\osp_{M|2n}$ and $\B_m(\om)$ (with $\om=M-2n$) on the tensor product space
$(\CC^{M|2n})^{\ot m}$ commute, thus allowing us to apply a similar approach to the
orthosymplectic Lie superalgebra.
However, the arguments used for the Lie algebras
do not readily extend to the super case. The main obstacle
is the singularity of the symmetriser $s^{(m)}$; it is a rational function of $\om$
which may have a pole
at $\om=M-2n$.

We solve the problem by using a new Brauer-type algebra $\BB_{2m+1}(\om)$ containing
$\B_{2m+1}(\om)$ as a subalgebra. We construct {\em abstract Segal--Sugawara vectors} as
elements of $\BB_{2m+1}(\om)$, keeping $\om$ as an indeterminate.
Then we show that the abstract Segal--Sugawara vectors admit an equivalent `integral form' where
the symmetriser $s^{(m)}$ is replaced by the symmetriser in the symmetric group algebra,
thus allowing the required evaluation of $\om$.

Our main result is Theorem~\ref{thm:osp}; it
provides explicit formulas for Segal--Sugawara vectors $\Phi_m$ for
$\osp_{M|2n}$. In particular, these elements are even, and they
generate a commutative subalgebra
of the universal enveloping algebra of
the Lie superalgebra $t^{-1}\osp_{M|2n}[t^{-1}]$; see Corollary~\ref{cor:commute}.
We believe that if $M$ is odd, then the vectors $\Phi_{2k}$
generate the centre of the affine vertex algebra,
as stated in Conjecture~\ref{con:gener}.

With the significance of the Feigin--Frenkel centre in applications
to Gaudin models and Vinberg's quantisation problem \cite{ffr:gm}, \cite{fft:gm}, \cite{r:si},
we expect that the explicit Segal--Sugawara vectors $\Phi_m$
will play a due role in understanding
the higher Hamiltonians in the Gaudin models with the $\osp$-symmetry and
the quantum shift of argument subalgebras in $\U(\osp_{M|2n})$.

\section{Segal--Sugawara vectors}
\label{sec:ssv}

We use the
involution $i\mapsto i\pr=M+2n-i+1$ on
the set $\{1,2,\dots,M+2n\}$. Set
\ben
\bi=\begin{cases} 1\qquad\text{for}\quad i=1,\dots,n,n',\dots,1',\\
0\qquad\text{for}\quad i=n+1,\dots,(n+1)\pr
\end{cases}
\een
and
\ben
\ta_i=\begin{cases} \phantom{-}1\qquad\text{for}\quad i=1,\dots,M+n,\\
-1\qquad\text{for}\quad i=M+n+1,\dots,M+2n.
\end{cases}
\een

\subsection{Lie superalgebras}
\label{subsec:ls}

A standard basis of the general linear Lie superalgebra $\gl_{M|2n}$ is formed by elements $E_{ij}$
of the parity $\bi+\bj\mod 2$ for $1\leqslant i,j\leqslant M+2n$, with the commutation relations
\ben
[E_{ij},E_{kl}]
=\de_{kj}\ts E_{i\tss l}-\de_{i\tss l}(-1)^{(\bi+\bj)(\bk+\bl)}\ts E_{kj}.
\een
We will regard the orthosymplectic Lie superalgebra $\osp_{M|2n}$
as the subalgebra
of $\gl_{M|2n}$ spanned by the elements
\ben
F_{ij}=E_{ij}-(-1)^{\bi\tss\bj+\bj}\ts\ta_i\ta_j E_{j'i'}.
\een
The corresponding affine Kac--Moody superalgebra $\widehat{\osp}_{M|2n} \coloneqq \osp_{M|2n}\left[t,t^{-1}\right] \oplus \C K$ has Lie superbracket
\begin{align}
\big[ F_{ij}[r], F_{kl}[s]\big]
    &{}= \delta_{jk}F_{il}[r+s]
- \delta_{il}(-1)^{(\bi+\bj)(\bk+\bl)}F_{kj} [r+s]
-  \delta_{ik'}(-1)^{\bi\bj + \bj} \varepsilon_i\varepsilon_j F_{j'l} [r+s] \nonumber \\[0.3em]
&{}+ \delta_{jl'}(-1)^{\bi\bar{k} +\bj\bar{k}}\varepsilon_i\varepsilon_j  F_{ki'}[r+s]
+r \delta_{r,-s} \left((-1)^{\bi} \delta_{il} \delta_{jk} - (-1)^{\bi\bj}\varepsilon_i \varepsilon_j \delta_{ik'} \delta_{jl'} \right) K,
\non
\end{align}
where $F_{ij}[r] \coloneqq F_{ij}t^r$ for each $r \in \Z$ and $1 \leq i,j \leq M+2n$, and $K$ is central.

\subsection{Vacuum module}

The vacuum module $V_{-h^{\vee}}(\osp_{M|2n})$ over $\widehat{\osp}_{M|2n}$
at the critical level is the quotient of the universal enveloping
algebra $\U(\widehat{\osp}_{M|2n})$ by the left ideal
generated by $K+h^\vee$ and $\osp_{M|2n}[t]$,
where $h^\vee=M-2n-2$ is the dual Coxeter number for $\osp_{M|2n}$.

Equipping the vacuum module $V_{-h^\vee}(\osp_{M|2n})$ with the derivation $\tau \coloneqq -d/dt$ yields a vertex algebra called the affine vertex algebra; see e.g. \cite{k:va} for details. As vector spaces, we have the isomorphism
\beql{isom}
V_{-h^\vee}(\osp_{M|2n}) \cong \U\big(t^{-1}\osp_{M|2n}[t^{-1}]\big).
\eeq
The centre $\z(\wh\osp_{M|2n})$ of the affine vertex algebra defined as in \eqref{ffcent} is called the Feigin--Frenkel centre, and its elements are called Segal--Sugawara vectors. Due to \eqref{isom}, the
Feigin--Frenkel centre can be regarded as a commutative subalgebra of
$\U\big(t^{-1}\osp_{M|2n}[t^{-1}]\big)$; see e.g. \cite[Sec.~3.3]{f:lc} and \cite[Sec.~6.2]{m:so}.
Moreover, this subalgebra is invariant under the derivation $\tau$ so that we can regard
$\z(\wh\osp_{M|2n})$ as a differential superalgebra.

\subsection{Main results}
\label{subsec:mr}

Consider the Lie superalgebra
$\widehat{\osp}_{M|2n}\oplus \C \tau$, where $\tau$ is even and satisfies
\ben
    [\tau, K]=0\Fand \big[\tau,F_{ij}[r]\big] = -rF_{ij}[r-1].
\een
We let $\U$ denote the universal enveloping algebra $\U\big(\widehat{\osp}_{M|2n}\oplus \C \tau\big)$.
Consider the tensor product superalgebra
\beql{tenprka}
\big(\End \C^{M|2n}  \big)^{\otimes m} \otimes \U,
\eeq
and for $1 \leq a \leq m$ and $r \in \Z$ define the elements
\beql{fra}
F[r]_a \coloneqq \sum_{i,j = 1}^{M+2n} 1 ^{\otimes (a-1)} \otimes e_{ij} \otimes 1^{\otimes (m-a)} \otimes F_{ij}[r] (-1)^{\bi \bj + \bi + \bj},
\eeq
where the $e_{ij}\in \End \C^{M|2n}$ are the standard matrix units.

The symmetric group $\Sym_m$ acts on the space $(\C^{M|2n})^{\ot m}$
by permuting tensor factors.
Denote by $H^{(m)}$ the element of the algebra
\eqref{tenprka} (with identity component in $\U$) which is
the image of the symmetriser $h^{(m)}\in\CC\Sym_m$
defined by
\beql{ha}
h^{(m)}=\frac{1}{m!}\sum_{s\in\Sym_m} s
\eeq
under the action of $\Sym_m$.
Furthermore,
let $\la=(\la_1,\dots,\la_{\ell})$ be a partition of $m$ of length $\ell=\ell(\la)$,
so that $\la_1\geqslant\dots\geqslant\la_{\ell}>0$ and $\la_1+\dots+\la_{\ell}=m$.
We denote by $c^{}_{\la}$ the number
of permutations in the symmetric group $\Sym_m$ of cycle type $\la$.
Set
\beql{fla}
F[-\la]=F[-\la_1]_1 \dots F[-\la_{\ell}]_{\ell},
\eeq
and for any $m\geq 2$, introduce elements $\Phi_m\in \U\big(t^{-1}\osp_{M|2n}[t^{-1}]\big)$ by
\beql{phimliesup}
\Phi^{}_{m}=\sum_{\la\ts\vdash m,\  \ell(\la)\ts \text{even}}\ts \Yc_{m,\ell}(M-2n-1)\ts
c^{}_{\la}\ts \str^{}_{1,\dots,\ell}\ts H^{(\ell)} F[-\la].
\eeq
Here we use the polynomials $\Yc_{m,\ell}(T)$ in a variable $T$ defined by
\ben
\Yc_{m,\ell}(T)=\frac{\ell!}{m!}\ts\prod_{k=\ell}^{m-1}(T+k);
\een
cf. \cite{r:iff}, while the supertrace
\ben
\str:\End \C^{M|2n}\to\CC,\qquad e_{ij}\mapsto \de_{ij}(-1)^{\bi}
\een
is taken over the first $\ell$ copies of $\End \C^{M|2n}$. The following is our main result.

\begin{theorem}\label{thm:osp}
All elements $\Phi_m$
belong to the Feigin--Frenkel centre $\z(\wh\osp_{M|2n})$.
\end{theorem}

The proof of the theorem will be given in Sec.~\ref{sec:as}, with some preliminary
results discussed in Sec.~\ref{sec:eb}. In the following corollary, we use the isomorphism
\eqref{isom}.

\begin{corollary}\label{cor:commute}
The elements $\Phi_m$ with $m=2,3,\dots$ generate a commutative subalgebra of $\U\big(t^{-1}\osp_{M|2n}[t^{-1}]\big)$.
\end{corollary}

\begin{conjecture}\label{con:gener}
If $M$ is odd, then the elements $\Phi^{}_{2}, \Phi^{}_{4}, \Phi^{}_{6},\dots$
generate $\z(\wh\osp_{M|2n})$ as a differential superalgebra.
\end{conjecture}

The conjecture holds for $n=0$ by \cite{m:ss}.
It is likely that for even values of
$M$ some additional elements of $\z(\wh\osp_{M|2n})$ arising from the super Pfaffian
are necessary to generate this superalgebra; cf. \cite{lz:io} and \cite{lw:sw}.

\begin{remark}\label{rem:mnze}
By taking $M=0$ or $n=0$ in Theorem~\ref{thm:osp},
we get elements of the respective Feigin--Frenkel centres
$\z(\wh\spa_{2n})$ or $\z(\wh\oa_{M})$. We thus obtain a new proof of the formulas
for the Segal--Sugawara vectors for the orthogonal and symplectic Lie algebras
given in \cite{m:ss}.
\qed
\end{remark}

\section{Extended Brauer-type algebra}
\label{sec:eb}

To prove Theorem~\ref{thm:osp}, we will consider its abstract version first.
Namely, we will introduce a new algebra $\BB_{2m+1}(\omega)$, and construct abstract analogues
of the Segal--Sugawara vectors $\phi_m$ as elements of $\BB_{2m+1}(\omega)$, regarding $\om$
as a variable. We show that the vectors satisfy the desired properties and that they are well-defined
at $\om=M-2n$. We then use a homomorphism from
$\BB_{2m+1}(M-2n)$ to the tensor product superalgebra to get the actual
Segal--Sugawara vectors $\Phi_m$ in the vacuum module over $\wh\osp_{M|2n}$.
The details will be given in Sec.~\ref{sec:as}, while this section
is devoted to some preliminary results on the algebra $\BB_{2m+1}(\omega)$.

\subsection{Brauer algebra}
\label{subsec:bra}

Let $\omega$ be an indeterminate. Recall that the Brauer algebra $\B_m(\omega)$ is the associative unital algebra with identity $1$, generated by $\su_a$ and $\eee_a$ with $1 \leq a \leq m-1$, subject only to the relations
\ben
\begin{aligned}
s_a^2&=1,\qquad \ep_a^2=\om\ts \ep_a,\qquad
\su_a\eee_a=\eee_a\su_a=\eee_a,\\[0.3em]
s_a s_b &=s_bs_a,\qquad \eee_a \eee_b = \eee_b \eee_a,\qquad
\su_a \eee_b = \eee_b\su_a,\qquad
|a-b|>1,\\[0.3em]
s_as_{a+1}s_a&=s_{a+1}s_as_{a+1},\qquad
\eee_a\eee_{a+1}\eee_a=\eee_a,\qquad \eee_{a+1}\eee_a\eee_{a+1}=\eee_{a+1},\\[0.3em]
\su_a\eee_{a+1}\eee_a&=\su_{a+1}\eee_a,\qquad \eee_{a+1}\eee_a
\su_{a+1}=\eee_{a+1}\su_a.
\end{aligned}
\een
For $1 \leq a<b \leq m$ we can also define
\begin{align}
\non
s_{ab}&=s_as_{a+1} \dots s_{b-2} s_{b-1} s_{b-2} \dots s_{a+1} s_a,\\
\epsilon_{ab} &=
s_a s_{a+1} \dots s_{b-2} \epsilon_{b-1} s_{b-2} \dots s_{a+1} s_a,
\non
\end{align}
and set $s_{ba}=s_{ab}$ and $\epsilon_{ba}=\epsilon_{ab}$.

The Brauer algebra $\B_m(\omega)$ also has a diagrammatic presentation as follows. The algebra has a basis of diagrams, where each diagram consists of two horizontal lines, each with $m$ nodes, and $m$ strings connecting the nodes pairwise. The product $xy$ of two diagrams $x$ and $y$ is computed by concatenation; we draw $y$ directly above $x$, connect the strings at the nodes in the middle, remove the middle line, and replace each loop formed by a factor of $\omega$. For example, in $\B_7(\omega)$, if
\begin{align}
    x\, &= \; \begin{tikzpicture}[baseline={([yshift=-1mm]current bounding box.center)},scale=0.55]
{
\draw (1,0) to[out=90,in=90] (3,0);
\draw (2,0) ..controls (2,1.2) and (7,1.2).. (7,0);
\draw (4,0) to[out=90,in=90] (6,0);
\draw (5,0) to[out=90,in=270] (6,2);
\draw (1,2) to[out=270,in=270] (4,2);
\draw (2,2) arc(180:360:0.5);
\draw (5,2) to[out=270,in=270] (7,2);
\draw[very thick] (0.5,0)--(7.5,0);
\draw[very thick] (0.5,2)--(7.5,2);
}
\end{tikzpicture}\; ,
&y \,&= \; \begin{tikzpicture}[baseline={([yshift=-1mm]current bounding box.center)},scale=0.55]
{
\draw (1,0) to[out=90,in=90] (4,0);
\draw (2,0) to[out=90,in=90] (5,0);
\draw (3,0) ..controls (3,0.9) and (7,0.9).. (7,0);
\draw (6,0) to[out=90,in=270] (4,2);
\draw (1,2) arc(180:360:0.5);
\draw (3,2) to[out=270,in=270] (5,2);
\draw (6,2) arc(180:360:0.5);
\draw[very thick] (0.5,0)--(7.5,0);
\draw[very thick] (0.5,2)--(7.5,2);
}
\end{tikzpicture}\; ,
\non
\end{align}
then
\begin{align}
    xy \, &= \;\;  \begin{tikzpicture}[baseline={([yshift=-1mm]current bounding box.center)},scale=0.55]
{
\begin{scope}[shift={(0,2)}]
\draw (1,0) to[out=90,in=90] (4,0);
\draw (2,0) to[out=90,in=90] (5,0);
\draw (3,0) ..controls (3,0.9) and (7,0.9).. (7,0);
\draw (6,0) to[out=90,in=270] (4,2);
\draw (1,2) arc(180:360:0.5);
\draw (3,2) to[out=270,in=270] (5,2);
\draw (6,2) arc(180:360:0.5);
\end{scope}
\draw (1,0) to[out=90,in=90] (3,0);
\draw (2,0) ..controls (2,1.2) and (7,1.2).. (7,0);
\draw (4,0) to[out=90,in=90] (6,0);
\draw (5,0) to[out=90,in=270] (6,2);
\draw (1,2) to[out=270,in=270] (4,2);
\draw (2,2) arc(180:360:0.5);
\draw (5,2) to[out=270,in=270] (7,2);
\draw[very thick] (0.5,0)--(7.5,0);
\draw[very thick] (0.5,2)--(7.5,2);
\draw[very thick] (0.5,4)--(7.5,4);
}
\end{tikzpicture}\;\;
= \; \omega^2 \; \begin{tikzpicture}[baseline={([yshift=-1mm]current bounding box.center)},scale=0.55]
{
\draw (1,2) arc(180:360:0.5);
\draw (3,2) to[out=270,in=270] (5,2);
\draw (6,2) arc(180:360:0.5);
\draw (1,0) to[out=90,in=90] (3,0);
\draw (2,0) ..controls (2,1.2) and (7,1.2).. (7,0);
\draw (4,0) to[out=90,in=90] (6,0);
\draw (5,0) to[out=90,in=270] (4,2);
\draw[very thick] (0.5,0)--(7.5,0);
\draw[very thick] (0.5,2)--(7.5,2);
}
\end{tikzpicture}\; .
\non
\end{align}
This can be identified with the non-diagrammatic presentation of $\B_m(\omega)$ by taking
\begin{align}
    s_{ab} &= \, \begin{tikzpicture}[baseline={([yshift=4.5mm]current bounding box.south)},scale=0.55]
{
\draw (0,0)--(0,2);
\draw (1,0)--(1,2);
\draw (3,0)--(3,2);
\draw (4,0)--(4,2);
\draw (6,0)--(6,2);
\draw (7,0)--(7,2);
\draw (2,0) to[out=90,in=270](5,2);
\draw (5,0) to[out=90,in=270](2,2);
\foreach \y in {0.2,0,-0.2}{
\node at (0.5+\y,1) {.};
}
\foreach \y in {0.2,0,-0.2}{
\node at (3.5+\y,0.5) {.};
}
\foreach \y in {0.2,0,-0.2}{
\node at (6.5+\y,1) {.};
}
\draw (0,2) node[font=\scriptsize,anchor=south]{$1$};
\draw (2,2) node[font=\scriptsize,anchor=south]{$a$};
\draw (5,2) node[font=\scriptsize,anchor=south]{$b$};
\draw (7,2) node[font=\scriptsize,anchor=south]{$m$};
\draw[very thick] (-0.5,0)--(7.5,0);
\draw[very thick] (-0.5,2)--(7.5,2);
}
\end{tikzpicture}\; ,
    &\epsilon_{ab} &= \, \begin{tikzpicture}[baseline={([yshift=4.5mm]current bounding box.south)},scale=0.55]
{
\draw (0,0)--(0,2);
\draw (1,0)--(1,2);
\draw (3,0)--(3,2);
\draw (4,0)--(4,2);
\draw (6,0)--(6,2);
\draw (7,0)--(7,2);
\draw (2,0) ..controls (2,0.9) and (5,0.9).. (5,0);
\draw (2,2) ..controls (2,1.1) and (5,1.1).. (5,2);
\foreach \y in {0.2,0,-0.2}{
\node at (0.5+\y,1) {.};
}
\foreach \y in {0.2,0,-0.2}{
\node at (3.5+\y,1) {.};
}
\foreach \y in {0.2,0,-0.2}{
\node at (6.5+\y,1) {.};
}
\draw (0,2) node[font=\scriptsize,anchor=south]{$1$};
\draw (2,2) node[font=\scriptsize,anchor=south]{$a$};
\draw (5,2) node[font=\scriptsize,anchor=south]{$b$};
\draw (7,2) node[font=\scriptsize,anchor=south]{$m$};
\draw[very thick] (-0.5,0)--(7.5,0);
\draw[very thick] (-0.5,2)--(7.5,2);
}
\end{tikzpicture}\; .
\non
\end{align}

For each $1 \leq a \leq m$, we define the \textit{partial transposition} $t_a:\B_m(\omega) \to \B_m(\omega)$ as the linear map taking each basis diagram $d$ to the diagram obtained by swapping the string endpoints connected to the nodes numbered $a$ on the top and bottom of $d$. For example, in $\B_6(\omega)$, we have
\begin{align}
    d &= \begin{tikzpicture}[baseline={([yshift=4.5mm]current bounding box.south)},scale=0.55]
{
\draw (0,0) to[out=90,in=270] (2,2);
\draw (1,0) to[out=90,in=270] (0,2);
\draw (2,0) ..controls (2,0.7) and (5,0.7).. (5,0);
\draw (3,0) to[out=90,in=270] (5,2);
\draw (4,0) to[out=90,in=270] (3,2);
\draw (1,2) to[out=270,in=270] (4,2);
\draw[very thick] (-0.5,0)--(5.5,0);
\draw[very thick] (-0.5,2)--(5.5,2);
}
\end{tikzpicture}\; ,
    &d\hspace{0.3mm}{}^{t_5} &= \begin{tikzpicture}[baseline={([yshift=4.5mm]current bounding box.south)},scale=0.55]
{
\draw (0,0) to[out=90,in=270] (2,2);
\draw (1,0) to[out=90,in=270] (0,2);
\draw (2,0) ..controls (2,0.7) and (5,0.7).. (5,0);
\draw (3,0) to[out=90,in=270] (5,2);
\draw (4,0) to[out=90,in=270] (1,2);
\draw (3,2) arc(180:360:0.5);
\draw[very thick] (-0.5,0)--(5.5,0);
\draw[very thick] (-0.5,2)--(5.5,2);
}
\end{tikzpicture}\; .
\non
\end{align}

\noindent
The {\em Brauer symmetriser} $s^{(m)} \in\B_{m}(\omega)$ is the unique nonzero element satisfying
\begin{align}
s^{(m)}s^{(m)} &= s^{(m)},
&s_as^{(m)} &= s^{(m)}s_a = s^{(m)},
&\epsilon_a s^{(m)} &= s^{(m)} \epsilon_a = 0
\non
\end{align}
for all $1 \leq a \leq m-1$. It follows that
\begin{align}
s_{ab}s^{(m)} &= s^{(m)} s_{ab} = s^{(m)},
&\epsilon_{ab}s^{(m)} &= s^{(m)} \epsilon_{ab} = 0
\non
\end{align}
for all $1 \leq a < b \leq m$. A few explicit expressions for $s^{(m)}$
are collected in \cite[Ch.~1]{m:so}. In particular,
\beql{symdia}
s^{(m)}=\frac{1}{m!}\ts\sum_{r=0}^{\lfloor m/2\rfloor}
(-1)^r\binom{\om/2+m-2}{r}^{-1}
\sum_{d\in\Dc^{(r)}}d,
\eeq
where $\Dc^{(r)}\subset\Bc_m(\om)$ denotes the set of diagrams
which have exactly $r$ horizontal strings in the top row; see
\cite{hx:ts}.

The Brauer algebra $\B_{m}(\omega)$ has a subalgebra isomorphic to the symmetric group algebra $\C\mathfrak{S}_m$, generated by $s_a$ with $1 \leq a \leq m-1$.
We thus have elements $h^{(m)} \in \B_{m}(\omega)$
defined as in \eqref{ha}.

Let $J_m$ be the vector subspace of $\B_m(\omega)$ spanned by sums of the form $d + d^{t_a}$, for any basis diagram $d$ and $1 \leq a \leq m$. From Lemma 2.4 of \cite{m:ss}, the Brauer and symmetric group symmetrisers in $\B_m(\omega)$ satisfy
\beql{equivjm}
    \gamma_m(\omega) s^{(m)} \equiv h^{(m)} \mod J_m,
\eeq
where
\beql{gam}
    \gamma_m(\omega) \coloneqq \frac{\omega + m - 2}{\omega + 2m - 2}.
\eeq
This result is readily generalised using the embeddings $\B_k(\omega) \hookrightarrow \B_m(\omega)$ for $1 \leq k \leq m$ obtained by adding $m-k$ pairs of nodes on the right of each diagram, connected by vertical strings. That is, in $\B_m(\omega)$, we have
\begin{align}
    \gamma_k(\omega) s^{(k)} \equiv h^{(k)} \mod J^{(k)}_m, \label{eq:gammash}
\end{align}
where $J^{(k)}_m$ is the vector subspace of $\B_m(\omega)$ spanned by sums of the form $d+d^{t_a}$, for any basis diagram $d$ whose rightmost $m-k$ pairs of nodes are connected by vertical strings, and
$1 \leq a \leq k$.

The Brauer algebra $\B_{m}(\om)$ with parameter $\om=M-2n$ acts on the tensor product
space $(\C^{M|2n})^{\otimes m}$ so that there is a homomorphism
\beql{brahom}
\B_{m}(M-2n) \to \big(\End\C^{M|2n}\big)^{\otimes m}
\eeq
which is defined by $s_{ab} \mapsto P_{ab}$ and $\epsilon_{ab} \mapsto Q_{ab}$,
where
\begin{align}
P_{ab} &\coloneqq \sum_{i,j=1}^{M+2n} 1^{\otimes (a-1)} \otimes e_{ij} \otimes 1^{\otimes (b-a-1)} \otimes e_{ji} \otimes 1^{\otimes (m-b)} (-1)^{\bj},
\non\\
Q_{ab} &\coloneqq \sum_{i,j=1}^{M+2n} 1^{\otimes (a-1)} \otimes e_{ij} \otimes 1^{\otimes (b-a-1)} \otimes e_{i'j'} \otimes 1^{\otimes (m-b)} (-1)^{\bi\bj + \bi + \bj} \varepsilon_i\varepsilon_j
\non
\end{align}
for $1 \leq a < b \leq m$. As we will work with the extended tensor product superalgebra,
we will usually identify these elements with
$P_{ab} \otimes 1$ and $Q_{ab} \otimes 1$ in \eqref{tenprka}, respectively.
For $b<a$, define $P_{ab} \coloneqq P_{ba}$ and $Q_{ab} \coloneqq Q_{ba}$.

The following easily verified property of the operators $Q_{ab}$ will be essential in the proof
of Theorem~\ref{thm:osp}. Suppose that $X$ is an
element of the superalgebra \eqref{tenprka}. We will also identify $X$ with the element $X\ot 1$
of this superalgebra but
with the parameter $m$ replaced by $m+1$.
Then for any $1\leq a\leq m$ we have
\beql{qxq}
Q_{a\ts m+1}X Q_{a\ts m+1}=(\str_a X)\ts Q_{a\ts m+1},
\eeq
where $\str_a$ denotes the partial supertrace taken over the $a$-th copy of the endomorphism superalgebra
in \eqref{tenprka}. We will use a version of relation \eqref{qxq} in the extended Brauer algebra
to define the `supertrace' on the algebra.

\subsection{Affine extension of the Brauer algebra}
\label{subsec:aff}

Our definition of the extended Brauer-type algebra is motivated by the matrix form
of the defining relations for the affine Kac--Moody superalgebra $\widehat{\osp}_{M|2n}$;
see Sec.~\ref{subsec:ls}. Namely, we can derive from the
defining relations of $\wh\osp_{M|2n}$ described in Sec.~\ref{subsec:ls}
that the elements $F[r]_a$ defined in \eqref{fra} satisfy the relations
\begin{multline}\non
F[r]_a F[s]_b- F[s]_b F[r]_a\\[0.4em]
{}= (P_{ab} - Q_{ab}) F[r+s]_b - F[r+s]_b (P_{ab}
- Q_{ab} ) +r\delta_{r,-s} (P_{ab} - Q_{ab} )K
\end{multline}
for all $1\leq a<b\leq m$.
Moreover, we also have $Q_{ab} F[r]_a Q_{ab}=0$,
\ben
P_{ab}F[r]_a=F[r]_b P_{ab}
\Fand Q_{ab} \big(F[r]_a+F[r]_b\big)=\big(F[r]_a+F[r]_b\big)Q_{ab}=0,
\een
whereas $P_{ab}$ and $Q_{ab}$ commute with $F[r]_c$ for $c\ne a,b$.

Following the approach of \cite{n:yo}, where
the degenerate affine Wenzl algebra was introduced (it is also known as the
Nazarov--Wenzl algebra), we use these matrix relations
to give the following definition.

\begin{definition}\label{def:affbra}
Let $\om$ be an indeterminate. Define
$\B^{\ts\aff}_m(\om)$ as
the associative unital algebra generated by the Brauer algebra $\B_m(\om)$
together with additional elements
$f[r]_a$ with $1\leq a\leq m$ and $r$ running over $\ZZ$, and $\tau, K$,
subject to the following relations. The element $K$ is central, $\tau$ commutes
with any element of $\B_m(\om)$, and we also have
\begin{multline}\label{frafsbabs}
f[r]_a f[s]_b- f[s]_b f[r]_a\\[0.4em]
{}= \left(s_{ab} - \epsilon_{ab}\right) f[r+s]_b - f[r+s]_b \left(s_{ab}
- \epsilon_{ab} \right) +r\delta_{r,-s} \left(s_{ab} - \epsilon_{ab} \right)K,
\end{multline}
\ben
\bal
s_{ab}f[r]_a=f[r]_b s_{ab}&\Fand \ep_{ab} \big(f[r]_a+f[r]_b\big)
=\big(f[r]_a+f[r]_b\big)\ep_{ab}=0,\\[0.4em]
s_{ab}f[r]_c=f[r]_c s_{ab}&\Fand \ep_{ab} f[r]_c=f[r]_c\ep_{ab} \qquad\text{for $c\ne a,b$},
\eal
\een
\medskip
\beql{taufco}
\epsilon_{ab} f[r]_a \epsilon_{ab} = 0,\qquad
f[r]_a \tau - \tau f[r]_a= r f[r-1]_a.
\eeq
\end{definition}

We have the epimorphism $\B^{\ts\aff}_m(\om)\to \B_m(\om)$ identical on elements of $\B_m(\om)$ and
sending $f[r]_a,\tau$ and $K$ to zero. Hence, we have a natural embedding
$\B_m(\om)\hra \B^{\ts\aff}_m(\om)$ so that $\B_m(\om)$  can be regarded as a subalgebra of
$\B^{\ts\aff}_m(\om)$.

The defining relations of $\B^{\ts\aff}_m(\om)$ show that this algebra is defined
over the algebra of polynomials $\CC[\om]$. However, we will also need to consider it over
the field of rational functions $\CC(\om)$,
as occurs already in \eqref{symdia}. We will not introduce
new notation for the extended algebra as this will be specified in the context.

\begin{remark}\label{rem:nw}
As with the Nazarov--Wenzl algebra, it should be natural to impose some additional relations in
$\B^{\ts\aff}_m(\om)$ to develop its reasonable structure theory; cf. \cite[Sec.~4]{n:yo}.
We will leave this topic outside the current paper,
as this algebra will play only an auxiliary role here.
\qed
\end{remark}

We will need the following key property of $\B^{\ts\aff}_m(\om)$.
Recall the tensor product superalgebra defined in \eqref{tenprka}.

\begin{proposition}\label{prop:hom}
The map that sends $\tau$ and $K$ to the elements with the same names, and
sends
\ben
s_{ab}\mapsto P_{ab},\qquad \ep_{ab}\mapsto Q_{ab},\qquad f[r]_a\mapsto F[r]_a,
\een
defines a homomorphism
\ben
\B^{\ts\aff}_m(M-2n)\to \big(\End \C^{M|2n}  \big)^{\otimes m} \otimes \U.
\een
\end{proposition}

\bpf
The verification of the relations is straightforward; they follow from the matrix form
of the defining relations in $\U$ pointed out above.
\epf

\subsection{Cyclic properties in the extended Brauer algebra}
\label{subsec:absss}

We will work with a slightly modified version of the algebra $\B^{\ts\aff}_m(\om)$.
Define $\BB_{2m+1}(\omega)$ as the associative unital algebra, generated by
the Brauer algebra $\B_{2m+1}(\om)$ with the elements $s_{ab}$ and $\ep_{ab}$
labelled by $a,b\in\{0,1,\dots,2m\}$, $a\ne b$, together with the
additional elements
$f[r]_a$ with $0\leq a\leq m$ and $r$ running over $\ZZ$, and $\tau, K$,
subject to the same relations as in Definition~\ref{def:affbra}.
The formulas as in Proposition~\ref{prop:hom} define
a homomorphism
\beql{homrho}
\rho: \BB_{2m+1}(M-2n) \to \left(\End\C^{M|2n}\right)^{\otimes (2m+1)} \otimes \U,
\eeq
where the copies of the endomorphism algebra are labelled by $0,1,\dots,2m$.
In $\BB_{2m+1}(\omega)$, for $1 \leq k \leq m$, set
\beql{qk}
\q{k} = \epsilon_{1\, m+1} \dots \epsilon_{k\, m+k}.
\eeq
We wish to use left- and right-multiplication by $\q{k}$ as an analogue of the partial trace operator
$\str_{1,\dots,k}$; cf. \eqref{qxq}.
We will need two cyclic properties in $\BB_{2m+1}(\om)$ described
in the following propositions. In the notation
below, the angle brackets indicate the subalgebras generated by the listed elements.

\begin{proposition}[Cyclic property 1]\label{prop:cyclic1}
For $1\leq k \leq m$, suppose
\begin{align}
    x &\in \e{s_a,\epsilon_a \mid 1 \leq a \leq k-1} =\B_k(\om)\subseteq \BB_{2m+1}(\omega),
    \non\\[0.4em]
    y &\in \e{s_a, \epsilon_a, f[r]_b, \tau, K \mid 0 \leq a \leq k-1,\ 0\leq b \leq k,\ r \in \Z} \subseteq \BB_{2m+1}(\omega).
    \non
\end{align}
Then
\ben
\q{k} xy  \q{k}
=
\q{k} yx \q{k}.
\een
\end{proposition}

\begin{proof}
We first note that, if we wish to verify a relation in $\BB_{2m+1}(\omega)$ that involves only the Brauer generators, it suffices to verify the corresponding relation in $\B_{2m+1}(\omega)$, because of the embedding $\B_{2m+1}(\omega) \hra \BB_{2m+1}(\om)$; see Sec.~\ref{subsec:aff}.
In particular, we can make use of the diagrammatic presentation of $\B_{2m+1}(\om)$.

For $1 \leq a < b \leq m$, then, it is easy to verify that
\begin{align}
\epsilon_{a\, m+a}\epsilon_{b\, m+b} s_{ab}
&= \epsilon_{a\, m+a} \epsilon_{b\, m+b} s_{m+a\, m+b} , \label{eq:eamaebmbsab}\\
\epsilon_{a\, m+a}\epsilon_{b\, m+b} \epsilon_{ab}
&= \epsilon_{a\, m+a} \epsilon_{b\, m+b} \epsilon_{m+a\, m+b}, \label{eq:eamaebmbeab}\\
s_{m+a\, m+b} \epsilon_{a\, m+a} \epsilon_{b\, m+b}
&= s_{ab} \epsilon_{a\, m+a} \epsilon_{b\, m+b}, \label{eq:sabeamaebmb}\\
\epsilon_{m+a\, m+b} \epsilon_{a\, m+a} \epsilon_{b\, m+b}
&= \epsilon_{ab} \epsilon_{a\, m+a} \epsilon_{b\, m+b}.\label{eq:eabeamaebmb}
\end{align}

Now, note that all of the indices of the $\epsilon$ generators involved in $\q{k}$ are distinct, so those generators commute. It follows that the first generator of $x$, either $s_{ab}$ or $\epsilon_{ab}$, can be converted to $s_{m+a\, m+b}$ or $\epsilon_{m+a\, m+b}$, using \eqref{eq:eamaebmbsab} or \eqref{eq:eamaebmbeab}, respectively. The resulting generator commutes with the rest of $x$ and all of $y$, and can be converted back to the generator we started with using \eqref{eq:sabeamaebmb} or \eqref{eq:eabeamaebmb}, respectively. Thus the first generator of $x$ can be moved to the right of $y$. This process can be repeated for all generators of $x$, because each $s_{ab}$ or $\epsilon_{ab}$ becomes $s_{m+a\, m+b}$ or $\epsilon_{m+a\, m+b}$ respectively, and then commutes with all other generators of $x$, and with $y$. The generators of $x$ can thus be moved to the right of $y$, ending up in the same order as they were in $x$, which proves the result.
\end{proof}

\begin{proposition}[Cyclic property 2] \label{prop:cyclic2}
For $1 \leq k \leq m$, suppose
\begin{align}
    x &\in \e{s_a, \epsilon_a, f[r]_b, \tau, K \mid 1 \leq a \leq k-1,\ 1\leq b \leq k,\ r \in \Z} \subseteq \BB_{2m+1}(\omega).
    \non
\end{align}
Then for any $1 \leq a \leq k$, and $y \in \{ s_{0a},\epsilon_{0a}\}$, we have
\begin{align*}
    \q{k} xy \q{k} = \q{k} yx \q{k}.
\end{align*}
\end{proposition}

\begin{proof}
One can check diagrammatically that, for $1 \leq a \leq k$,
\begin{alignat}{2}
\epsilon_{a\, m+a} s^{}_{0a} &= \epsilon_{a\, m+a} \epsilon^{}_{0\, m+a},
&\epsilon^{}_{0\, m+a} \epsilon_{a\, m+a} &= s^{}_{0a} \epsilon_{a\, m+a},
\non\\
\epsilon_{a\, m+a} \epsilon^{}_{0a} &= \epsilon_{a\, m+a} s^{}_{0\, m+a},
\qquad\qquad
&s^{}_{0\, m+a} \epsilon_{a\, m+a} &= \epsilon^{}_{0a} \epsilon_{a\, m+a}.
\non
\end{alignat}
It follows that
\ben
\q{k}s^{}_{0a}x \q{k} = \q{k} \epsilon^{}_{0\, m+a} x \q{k}
= \q{k} x \epsilon^{}_{0\, m+a} \q{k}
= \q{k} x s^{}_{0a} \q{k},
\een
since each generator of $x$ commutes with $\epsilon^{}_{0\, m+a}$. The relation for $y=\epsilon_{0a}$ is proven analogously.
\end{proof}

\section{Abstract Segal--Sugawara vectors}
\label{sec:as}

Now we explain our strategy to prove
Theorem~\ref{thm:osp} in more detail. The theorem will follow by verifying that
the elements $\Phi_m\in V_{-h^\vee}(\osp_{M|2n})$ are annihilated by the action
of $\osp_{M|2n}[t]$; see \eqref{ffcent}. Since the orthosymplectic Lie superalgebra
is simple, it suffices to show that
\beql{annihot}
F_{ij}[0]\Phi_m=F_{ij}[1]\Phi_m=0
\eeq
for all $i,j$. In the case $n=0$ (i.e., for the orthogonal Lie algebra $\oa_M$),
the verification of \eqref{annihot} was carried out
with the use of the matrix techniques, where $\Phi_m$ is first represented as a weighted trace
of the symmetriser in the Brauer algebra \cite{m:ff}, and then brought
to the form \eqref{phimliesup}; see \cite{m:ss}.
In the super case for arbitrary $n$, the corresponding initial expression would take the form
\beql{defphi}
\Phi^{}_{m}=\ga_m(M-2n)\ts\str^{}_{1,\dots,m}\ts S^{(m)}
\big(\tau+F[-1]_1\big)\dots \big(\tau+F[-1]_m\big)\tss 1,
\eeq
where $\ga_m(\om)$ is defined in \eqref{gam}, while $S^{(m)}$ denotes the image of the symmetriser
$s^{(m)}$ introduced in \eqref{symdia}, under the
homomorphism \eqref{brahom}, and we assume $\tau 1=0$. However, since
$\ga_m(\om)$ and $s^{(m)}$ are rational functions in $\om$,
the expression in \eqref{defphi}
is not defined for some values of $M$ and $n$.
The same problem already occurs for the symplectic Lie algebra $\spa_{2n}$ (i.e., for $M=0$);
a way around this was found in \cite{m:ff} (see also \cite[Sec.~8.3]{m:so})
with the use of `analytic continuation' over $n$.
Its extension to arbitrary $M$ appears to be problematic as
a super version of
the argument should rely on
orthosymplectic invariant theory, which is much less well-understood; cf. \cite{lz:io}.

Our way to settle the singularity issue of the expression in \eqref{defphi}
is to `lift' it to the algebra
$\BB_{2m+1}(\om)$ while keeping $\om$ as an indeterminate, thus working over the field
of rational functions $\CC(\om)$. Namely, we find some
`abstract' counterparts $\phi_m\in\BB_{2m+1}(\om)$ of $\Phi_m$
satisfying the desired
annihilation properties. This requires consistent definitions of the supertraces
in $\BB_{2m+1}(\om)$ and in the tensor product superalgebra appearing in \eqref{homrho}.
For the latter, we use the observation
\eqref{qxq}, which implies a counterpart of \eqref{defphi} without an explicit use
of the supertrace:
\ben
\Phi_m\Q{m}=\Q{m}\big(\ga_m(M-2n)\ts S^{(m)}
\big(\tau+F[-1]_1\big)\dots \big(\tau+F[-1]_m\big)\tss 1\ts\big) \Q{m},
\een
assuming $\Phi_m$ is defined,
where
for each $1 \leq k \leq m$ we set
\begin{align}
    \Q{k} &\coloneqq Q_{1\, m+1} Q_{2\, m+2} \dots Q_{k\, m+k} \in \left(\End \C^{M|2n}\right)^{\otimes 2m}.
    \non
\end{align}
As a final step,
we use \eqref{equivjm}
to find an equivalent `integral form' of the abstract Segal--Sugawara vectors $\phi_m$,
implying that the singularities in \eqref{defphi} are removable, by
showing that the vectors $\phi_m$
allow for a well-defined evaluation at $\om=M-2n$, yielding formula \eqref{phimliesup}.

\subsection{Annihilation properties}
\label{subsec:ap}

To implement the program,
begin by setting
\begin{align}
\f_a \coloneqq \tau + f[-1]_a \in \BB_{2m+1}(\omega),\qquad a=1,\dots,m.
\non
\end{align}
Use notation \eqref{qk} and introduce the
element
\beql{absrat}
\q{m} s^{(m)} \f_1 \dots \f_m \q{m}\in \BB_{2m+1}(\omega).
\eeq
We can regard it as a polynomial in $\tau$ by moving
the powers of $\tau$ to their right-most positions by using
the second relation in \eqref{taufco}.
Our goal is to show that
all coefficients of this polynomial are
{\em abstract Segal--Sugawara vectors}
in the sense that both expressions
\beql{eq:f1sexpression}
f[0]_0 \ts\q{m} s^{(m)} \f_1 \dots \f_m \q{m} \Fand
f[1]_0 \ts\q{m} s^{(m)} \f_1 \dots \f_m \q{m}
\eeq
are zero modulo the left ideal
in $\BB_{2m+1}(\omega)$ generated by the subspaces $f[r]_a\tss \B_{2m+1}(\om)$ for $r\geq 0$
and $a=0,\dots,m$, assuming that the `abstract level' is critical: $K=-\om+2$.

\begin{proposition}\label{prop:thm816first}
In $\BB_{2m+1}(\omega)$, we have
\ben
f[0]_0 \ts\q{m} s^{(m)} \f_1 \dots \f_m \q{m}= \q{m} s^{(m)} \f_1 \dots \f_m \q{m} f[0]_0.
\een
\end{proposition}

\begin{proof}
Let $\vp_{xy} \coloneqq s_{xy} - \epsilon_{xy}$. For $1 \leq a \leq m$, we have from
the defining relations in $\BB_{2m+1}(\om)$ and \eqref{frafsbabs} that
\ben
f[0]_0\f_a - \f_a f[0]_0
= \vp_{0a} \f_a - \f_a \vp_{0a}.
\een
Then
\ben
\bal
f[0]_0 \q{m} s^{(m)} \f_1 \dots \f_m \q{m}
&{}=  \sum_{a=1}^m \q{m} s^{(m)}
\f_1 \dots \f_{a-1}\big(f[0]_0 \f_a - \f_a f[0]_0\big) \f_{a+1} \dots \f_m \q{m}\\
&{}\qquad + \q{m} s^{(m)} \f_1 \dots \f_m f[0]_0\q{m}.
\eal
\een
The sum equals
\begin{multline}
\sum_{a=1}^m \q{m} s^{(m)} \f_1 \dots \f_{a-1}
\big(\vp_{0a} \f_a - \f_a \vp_{0a}\big) \f_{a+1} \dots \f_m \q{m}\\
=  \sum_{a=1}^m \q{m} s^{(m)} \vp_{0a} \f_1 \dots \f_m \q{m}  - \sum_{a=1}^m \q{m} s^{(m)} \f_1 \dots   \f_m \vp_{0a}\q{m}. \nonumber
\end{multline}
Applying cyclic property 2, and noting
that $s^{(m)}$ commutes with the sum of $\vp_{0a}$, this becomes
\ben
\sum_{a=1}^M \q{m} \vp_{0a} s^{(m)} \f_1 \dots \f_m \q{m}  - \sum_{a=1}^M \q{m}  \vp_{0a} s^{(m)} \f_1 \dots   \f_m\q{m} =0,
\een
as desired.
\end{proof}

The second expression in \eqref{eq:f1sexpression} is more difficult to handle, requiring a few lemmas.

\begin{lemma}\label{lem:817}
In $\BB_{2m+1}(\omega)$, for $1 \leq k \leq m$, we have
\beql{eq:lem817result}
\epsilon_{k\, m+k} s^{(k)} \vp_{0k} \epsilon_{k\, m+k}
= \frac{\omega + 2k - 2}{k(\omega + 2k - 4)} s^{(k-1)} \bigg(\sum_{a=1}^{k-1} \vp_{0a}\bigg)
\epsilon_{k\, m+k}.
\eeq
\end{lemma}

\begin{proof}
In the Brauer algebra $\B_k(\omega)$, it is known that
\ben
s^{(k)} = \frac{1}{k(\omega + 2k -4)} \bigg(1 + \sum_{a=1}^{k-1} (s_{ak} - \epsilon_{ak}) \bigg)
\bigg(\omega + k-3 + \sum_{a=1}^{k-1}(s_{ak} - \epsilon_{ak}) \bigg) s^{(k-1)};
\een
see \cite[proof of Lemma 1.3.2]{m:so}. Expanding this and applying Brauer relations, we find
\beql{eq:symrecursion2}
s^{(k)}
= \frac{1}{k} \bigg(1 +
\sum_{a=1}^{k-1} s_{ak} -
\frac{2}{\omega + 2k - 4} \bigg(\sum_{a=1}^{k-1} \epsilon_{ak}  + \sum_{1 \leq a < b \leq k-1}   s_{ak}\epsilon_{bk}\, \bigg)
 \!\bigg)s^{(k-1)}.
\eeq
Due to the embeddings $\B_k(\omega) \hra \B_m(\omega) \hra \BB_{2m+1}(\omega)$,
this result holds in $\BB_{2m+1}(\omega)$ also.

Substituting the above expression for $s^{(k)}$ into the left-hand side of \eqref{eq:lem817result}, we can commute the $\vp_{0k}$ left past the $\epsilon_{k\, m+k}$ and $s^{(k-1)}$, expand, and simplify the resulting expression by noting that
\begin{alignat}{2}
\epsilon_{k\, m+k}\ts \vp_{0k}\ts  \epsilon_{k\, m+k} &= 0, &
\epsilon_{k\, m+k}\ts  s_{ak}\ts  \vp_{0k}\ts  \epsilon_{k\, m+k} &= \vp_{0a}\ts  \epsilon_{k\, m+k},
\non\\
\epsilon_{k\, m+k}\ts  \epsilon_{ak}\ts  \vp_{0k}\ts  \epsilon_{k\, m+k} &=
- \vp_{0a}\ts  \epsilon_{k\, m+k},
&\qquad\epsilon_{k\, m+k}\ts  s_{ak}\ts  \epsilon_{bk}\ts  \vp_{0k}\ts  \epsilon_{k\, m+k} &= \epsilon_{ab} \ts \ts \vp_{0a}\ts
\epsilon_{k\, m+k}.
\non
\end{alignat}
This gives
\begin{align}
\epsilon_{k\, m+k}\ts  s^{(k)}\ts  \vp_{0k}\ts  \epsilon_{k\, m+k}
&= \frac{1}{k}\ts \bigg(\frac{\omega + 2k - 2}{\omega + 2k - 4} \,\sum_{a=1}^{k-1} \vp_{0a}
-
\frac{2}{\omega + 2k - 4} \sum_{1 \leq a < b \leq k-1}   \epsilon_{ab} \vp_{0a}
 \bigg) \epsilon_{k\, m+k} s^{(k-1)}
 \non
 \end{align}
which equals
\ben
\frac{\omega + 2k - 2}{k(\omega + 2k - 4)}
\bigg(\sum_{a=1}^{k-1} \vp_{0a} \bigg)\epsilon_{k\, m+k} s^{(k-1)},
\een
since $s^{(k-1)}$ commutes with $\epsilon_{k\, m+k}$, and
\begin{align}
\epsilon_{ab}\ts  \vp_{0a}\ts  \epsilon_{k\, m+k}\ts  s^{(k-1)}
&= \epsilon_{ab}\ts  \vp_{0a}\ts  s^{(k-1)}\ts  \epsilon_{k\, m+k}=0.
\non
\end{align}
The latter holds by noting that
\ben
\epsilon_{ab}\ts  \vp_{0a}\ts  s^{(k-1)}
=\epsilon_{ab}\ts  \vp_{0a}\ts  s_{ab}s^{(k-1)}
=\epsilon_{ab}\ts  s_{ab}\ts \vp_{0b}\ts s^{(k-1)}=\epsilon_{ab}\ts \vp_{0b}\ts s^{(k-1)}=
-\epsilon_{ab}\ts  \vp_{0a}\ts  s^{(k-1)},
\een
completing the proof.
\end{proof}

\begin{lemma}\label{lem:lemma132}
In $\BB_{2m+1}(\omega)$, for $1 \leq k \leq m$ we have
\begin{align}
\epsilon_{k\, m+k}\ts s^{(k)} \epsilon_{k\, m+k} = \frac{(\omega + k - 3)(\omega + 2k -2)}{k(\omega + 2k - 4)} s^{(k-1)} \epsilon_{k\, m+k}.
\non
\end{align}
\end{lemma}

\begin{proof}
We substitute the expression \eqref{eq:symrecursion2} for $s^{(k)}$ into the left-hand side, and
then use that $s^{(k-1)}$ and $\epsilon_{k\, m+k}$ commute to pull both copies of $\epsilon_{k\, m+k}$ inside the brackets. We then simplify each term using the Brauer relations, and pull a remaining copy of $\epsilon_{k\, m+k}$ back out of the brackets, to the right. This gives
\begin{align}
&\epsilon_{k\, m+k}\ts s^{(k)} \epsilon_{k\, m+k} \nonumber \\
&= \frac{1}{k} \Big(\omega  +
\sum_{a=1}^{k-1} 1 -
\frac{2}{\omega + 2k - 4} \sum_{a=1}^{k-1} 1  -
\frac{2}{\omega + 2k - 4}\sum_{1 \leq a < b \leq k-1}   \epsilon_{ab} \Big)s^{(k-1)} \epsilon_{k\, m+k},
\non
\end{align}
which may be simplified to give the result by noting that $\epsilon_{ab} s^{(k-1)} = 0$.
\end{proof}

\begin{lemma}\label{lem:Maninconsequence}
Let $1 \leq a < b \leq k \leq m$. Then, in $\BB_{2m+1}(\omega)$, we have
\begin{align}
s_{ab}\ts \f_1 \dots \f_k s^{(k)}
= \f_1 \dots \f_k s^{(k)}.
\non
\end{align}
\end{lemma}
\begin{proof}
It suffices to show this for $b=a+1$, since $s_{ab}$ is equal to a product of generators $s_c$ with $a \leq c < b$. Using $\f_a = \tau + f[-1]_a$ and the defining relations of $\BB_{2m+1}(\omega)$, it is tedious but straightforward to show that
\begin{align}
\f_a \f_{a+1} \left(\frac{1+s_a}{2} - \frac{\epsilon_a}{\omega} \right)
= \f_{a+1} \f_a \left(\frac{1+s_a}{2} - \frac{\epsilon_a}{\omega} \right) .
\non
\end{align}
We note also that the expression in the brackets is absorbed by $s^{(k)}$, and commutes with $\f_b$ for $b\geq a+2$, while $s_a$ commutes with $\f_b$ for $b\leq a-1$. The result follows readily from these observations.
\end{proof}

\begin{corollary}\label{cor:lemma818}
Let $1 \leq k \leq m$. Then for any $1 \leq a < b \leq k$, we have
\begin{align}
\q{k} s^{(k)} \vp_{0a} \f_1 \dots \f_k \q{k}
&= \q{k} s^{(k)} \vp_{0b} \f_1 \dots \f_k \q{k}
\non
\end{align}
in $\BB_{2m+1}(\omega)$.
\end{corollary}
\begin{proof}
This follows from $s^{(k)}s_{ab} = s_{ab}\tss s^{(k)} = s^{(k)}$, Lemma \ref{lem:Maninconsequence}, and cyclic property 1.
\end{proof}

We are now in a position to establish the desired property of the second expression in
\eqref{eq:f1sexpression}, as given in the next proposition.

\begin{proposition} \label{prop:thm816second}
In $\BB_{2m+1}(\omega)$, we have
\begin{align}
&
f[1]_0\ts \q{m} s^{(m)} \f_1 \dots \f_m \q{m}
\nonumber \\[0.4em]
&= (\omega + K - 2)\ts\frac{\omega + 2m -2}{\omega + 2m -4}\ts
\bigg(\sum_{a=1}^{m-1}\q{m-1} s^{(m-1)}  \vp_{0a} \f_1 \dots \f_{m-1} \q{m-1} \bigg) \epsilon_{m\, 2m}\nonumber \\[0.4em]
&\qquad - m \q{m} \vp_{0m}s^{(m)}  \f_1 \dots \f_{m-1} f[0]_m \q{m}
+ m\q{m} s^{(m)}  \f_1 \dots \f_{m-1} \q{m} f[0]_0 \nonumber \\[0.6em]
&\qquad + m \q{m} s^{(m)} \vp_{0m} \f_1 \dots \f_{m-1} f[0]_m \q{m}
 + \q{m} s^{(m)} \f_1 \dots \f_m \q{m} f[1]_0.
\non
\end{align}
\end{proposition}

\begin{proof}
We first note that
\begin{align}
f[1]_0\f_a - \f_a f[1]_0
=f[0]_0 + \vp_{0a} f[0]_a - f[0]_a \vp_{0a} + \vp_{0a}K, \label{eq:f1Facomm}
\end{align}
and
\begin{align}
f[0]_a\f_b - \f_b f[0]_a
&= \vp_{ab} \f_b - \f_b \vp_{ab}. \label{eq:faFbcomm}
\end{align}
Similar to the proof of Proposition \ref{prop:thm816first}, we begin by rewriting our expression as a telescoping sum. We then apply \eqref{eq:f1Facomm}, finding
\begin{align}
&f[1]_0\ts \q{m} s^{(m)} \f_1 \dots \f_m \q{m}
\non\\
&= \sum_{a=1}^m \q{m} s^{(m)} \f_1 \dots \f_{a-1} \big(f[0]_0 + \vp_{0a} f[0]_a - f[0]_a \vp_{0a} + \vp_{0a}K\big) \f_{a+1} \dots \f_m \q{m} \nonumber \\
&\qquad + \q{m} s^{(m)} \f_1 \dots \f_m \q{m} f[1]_0.
\non
\end{align}
Expanding the brackets gives four sums; we rewrite the first three
as telescoping sums in an additional index $b$, and apply \eqref{eq:faFbcomm} to represent
the sum over $a$ in the form
\begin{align}
&\sum_{1 \leq a < b \leq m} \q{m} s^{(m)} \f_1 \dots \f_{a-1}\f_{a+1} \dots \f_{b-1}
\big(\vp_{0b} \f_b - \f_b \vp_{0b}\big) \f_{b+1}  \dots \f_m \q{m} \nonumber \\
&\qquad + \sum_{a=1}^m\q{m} s^{(m)} \f_1 \dots \f_{a-1}\f_{a+1} \dots \f_m \q{m} f[0]_0 \nonumber \\
&\qquad + \sum_{1\leq a < b \leq m} \q{m} s^{(m)} \vp_{0a} \f_1 \dots \f_{a-1}  \f_{a+1} \dots \f_{b-1}\big(\vp_{ab} \f_b - \f_b \vp_{ab}\big) \f_{b+1}\dots\f_m \q{m} \nonumber \\
&\qquad + \sum_{a=1}^m \q{m} s^{(m)} \vp_{0a} \f_1 \dots \f_{a-1}  \f_{a+1} \dots\f_m f[0]_a \q{m} \nonumber \\
&\qquad - \sum_{1 \leq a < b \leq m} \q{m} s^{(m)} \f_1 \dots \f_{a-1}  \f_{a+1} \dots \f_{b-1}
\big(\vp_{ab} \f_b - \f_b \vp_{ab}\big) \f_{b+1} \dots \f_m \vp_{0a}\q{m} \nonumber \\
&\qquad - \sum_{a=1}^m\q{m} s^{(m)} \f_1 \dots \f_{a-1} \f_{a+1} \dots  \f_m f[0]_a \vp_{0a}\q{m}
\non\\
&\qquad + K \sum_{a=1}^m \q{m} s^{(m)} \vp_{0a} \f_1 \dots \f_{a-1} \f_{a+1} \dots \f_m \q{m}.
\non
\end{align}
We now manipulate the three sums over $a$ and $b$, and the final sum. Since $s^{(m)}s_{xy} = s^{(m)}$ for $1 \leq x < y \leq m$, we can insert $s_{m-1\, m} \dots s_{a\, a+1}$ after the $s^{(m)}$ in each sum. We move this product to the right in each word, using
\ben
    s_{wx}\vp_{yz} = \vp_{yz} s_{wx}, \qquad
    s_{wx}\vp_{xy} = \vp_{wy} s_{wx}, \qquad
    s_{wx} \vp_{wx} = \vp_{wx} s_{wx},
\een
and
\ben
    s_{wx}\f_y = \f_y s_{wx}, \qquad
    s_{wx}\f_x = \f_w s_{wx}
\een
for $w,x,y,z$ distinct, and then apply cyclic property 1 to bring it back to the left of $s^{(m)}$, where it is absorbed. We also move $\vp_{xy}$ throughout the word, using $\vp_{xy}\f_z = \f_z \vp_{xy}$ for distinct $x,y,z$, and cyclic property 1 or 2 as appropriate; that is, for $x,y \neq 0$, and $x=0$, respectively. Additionally, for $1 \leq x<y\leq m$, we use $\vp_{xy} s^{(m)} = s^{(m)} \vp_{xy} = s^{(m)}$. The resulting summands are then independent of one or more of the summation indices, and we simplify the sums accordingly.

In the first sum, for example, inserting the product $s_{m-1\, m} \dots s_{a\, a+1}$ gives
\begin{align}
&\sum_{1 \leq a < b \leq m} \q{m} s^{(m)}s_{m-1\, m} \dots s_{a\, a+1} \f_1 \dots \f_{a-1}\f_{a+1} \dots \f_{b-1}\ts \vp_{0b}\ts \f_b  \f_{b+1}  \dots \f_m \q{m} \nonumber \\
&\qquad -\sum_{1 \leq a < b \leq m} \q{m} s^{(m)} s_{m-1\, m} \dots s_{a\, a+1}\f_1 \dots \f_{a-1}\f_{a+1} \dots \f_{b-1}  \f_b\ts \vp_{0b}\ts \f_{b+1}  \dots \f_m \q{m},
\non
\end{align}
which equals
\begin{align}
&\sum_{1 \leq a < b \leq m} \q{m} s^{(m)} \f_1 \dots \f_{b-2} \ts\vp_{0\, b-1}
\ts  \f_{b-1}  \dots \f_{m-1} s_{m-1\, m} \dots s_{a\, a+1} \q{m} \nonumber \\
&\qquad -\sum_{1 \leq a < b \leq m} \q{m} s^{(m)} \f_1 \dots \f_{b-2}  \f_{b-1}
\ts \vp_{0\, b-1}\ts \f_{b}  \dots \f_{m-1} s_{m-1\, m} \dots s_{a\, a+1} \q{m}.
\non
\end{align}
We then apply cyclic property 1 and the absorption property
of the symmetriser, followed by $\vp_{xy}\f_z = \f_z \vp_{xy}$
for distinct $x,y,z$, to write the first sum as
\begin{align}
& \sum_{1 \leq a < b \leq m} \big(\q{m} s^{(m)} \f_1 \dots \f_{b-2} \ts\vp_{0\, b-1}\ts
\f_{b-1}  \dots \f_{m-1} \q{m} -
\q{m} s^{(m)} \f_1 \dots  \f_{b-1} \ts\vp_{0\, b-1}\ts  \f_{b}  \dots \f_{m-1} \q{m}\big),
\non\\
&\qquad = \sum_{1 \leq a < b \leq m} \big(\q{m} s^{(m)} \ts\vp_{0\, b-1}\ts  \f_1 \dots \f_{m-1} \q{m}
- \q{m} s^{(m)} \f_1 \dots \f_{m-1} \ts\vp_{0\, b-1}\ts  \q{m}\big).
\non
\end{align}
The summands do not depend on $a$, so may be simplified to
\begin{align}
& \sum_{a=1}^{m-1} a \q{m} s^{(m)} \ts\vp_{0a}\ts  \f_1 \dots \f_{m-1} \q{m}
-\sum_{a=1}^{m-1} a\q{m} s^{(m)} \f_1 \dots \f_{m-1} \ts\vp_{0a}\ts  \q{m}
\non\\
&= \sum_{a=1}^{m-1} a \q{m} s^{(m)} \ts\vp_{0a}\ts  \f_1 \dots \f_{m-1} \q{m}
-\sum_{a=1}^{m-1} a\q{m} \ts\vp_{0a}\ts  s^{(m)} \f_1 \dots \f_{m-1} \q{m},
\non
\end{align}
where the last line uses cyclic property 2.

Simplifying the other sums similarly, and grouping some like terms, we represent
the initial expression $f[1]_0\ts \q{m} s^{(m)} \f_1 \dots \f_m \q{m}$ as
\begin{align}
&\sum_{a=1}^{m-1} a\q{m} s^{(m)}  \vp_{0a} \f_1 \dots \f_{m-1} \q{m}
- \sum_{a=1}^{m-1} a\q{m} \ts\vp_{0a}\ts s^{(m)}  \f_1 \dots \f_{m-1} \q{m} \nonumber \\
&\qquad + \sum_{a=1}^{m-1} a\q{m} s^{(m)} \vp_{0m}\vp_{a\, m}  \f_1 \dots \f_{m-1} \q{m}
+\sum_{a=1}^{m-1} a \q{m} \vp_{a\tss m}\vp_{0m} s^{(m)} \f_1 \dots \f_{m-1} \q{m} \nonumber
\end{align}
plus the sum of the terms
\begin{align}
&m(K-m+1)\q{m}s^{(m)} \vp_{0m}  \f_1 \dots \f_{m-1}  \q{m}
 - m \q{m} \vp_{0m}s^{(m)}  \f_1 \dots \f_{m-1} f[0]_m \q{m} \nonumber\\[0.4em]
&+ m\q{m} s^{(m)}  \f_1 \dots \f_{m-1} \q{m} f[0]_0  + m \q{m} s^{(m)} \vp_{0m} \f_1 \dots \f_{m-1} f[0]_m \q{m} \nonumber \\[0.4em]
&\qquad + \q{m} s^{(m)} \f_1 \dots \f_m \q{m} f[1]_0 .
\non
\end{align}
Now, using the Brauer relations and the properties of the symmetriser, we have
\ben
s^{(m)} \vp_{0m} s_{am}
= s^{(m)} s_{am} \ts\vp_{0a}\ts
= s^{(m)} \vp_{0a}
\een
and
\ben
s^{(m)} \vp_{0m} \epsilon_{am}
= s^{(m)} s_{0m} \epsilon_{am} - s^{(m)} \epsilon_{0m} \epsilon_{am}
= s^{(m)} s_{0m} \epsilon_{am} - s^{(m)} s_{0m} \epsilon_{am}
= 0,
\een
so that
\ben
s^{(m)} \vp_{0m} \vp_{am} = s^{(m)} \vp_{0a},
\een
and analogously,
\ben
\vp_{am} \vp_{0m} s^{(m)} = \vp_{0a} s^{(m)}.
\een
The four sums over $a$ may thus be simplified to
\begin{multline}
2\sum_{a=1}^{m-1} a\q{m} s^{(m)}  \vp_{0a} \f_1 \dots \f_{m-1} \q{m}
= 2\sum_{a=1}^{m-1} a\q{m-1} \epsilon_{m\, 2m}s^{(m)}
\epsilon_{m\, 2m} \vp_{0a} \f_1 \dots \f_{m-1} \q{m-1}
\non\\[0.4em]
{}= \frac{2(\omega + m - 3)(\omega + 2m -2)}{m(\omega + 2m - 4)}\sum_{a=1}^{m-1} a\q{m-1}  s^{(m-1)}  \vp_{0a} \f_1 \dots \f_{m-1} \q{m-1} \epsilon_{m\, 2m},
\end{multline}
using $\q{m} = \q{m-1}\epsilon_{m\, 2m} = \epsilon_{m\, 2m} \q{m-1}$ and Lemma~\ref{lem:lemma132} with $k=m$. From Corollary~\ref{cor:lemma818} with $k=m-1$, however, $\q{m-1} s^{(m-1)}\vp_{0a}   \f_1 \dots \f_{m-1}  \q{m-1}$ has the same value for each $a$, so this can be further simplified to
\begin{align}
&\frac{(\omega + m - 3)(\omega + 2m -2)}{\omega + 2m - 4} \sum_{a=1}^{m-1}
\q{m-1}  s^{(m-1)}  \vp_{0a} \f_1 \dots \f_{m-1} \q{m-1} \epsilon_{m\, 2m}.
\non
\end{align}
Similarly, using Lemma~\ref{lem:817} with $k=m$, we also have
\begin{align}
&m(K-m+1)\q{m}s^{(m)} \vp_{0m}  \f_1 \dots \f_{m-1}  \q{m}
\non\\[0.4em]
&= \frac{(K-m+1)(\omega + 2m -2)}{\omega+2m-4}  \sum_{a=1}^{m-1} \q{m-1} s^{(m-1)}\vp_{0a}
\f_1 \dots \f_{m-1}  \q{m-1} \epsilon_{m\, 2m}.
\non
\end{align}
Substituting these results back into the expansion of $f[1]_0\ts \q{m} s^{(m)} \f_1 \dots \f_m \q{m}$
 gives the desired expression.
\end{proof}

\subsection{Integral form of the abstract Segal--Sugawara vectors}
\label{subsec:rf}

Propositions~\ref{prop:thm816first} and \ref{prop:thm816second} show
that all coefficients of the polynomial in $\tau$
defined in \eqref{absrat} are abstract Segal--Sugawara vectors. However, the coefficients
are defined in the algebra $\BB_{2m+1}(\omega)$ over the field $\CC(\om)$.
Since we aim to evaluate $\om$ at $M-2n$, we would like to bring the coefficients
to an `integral form' to be regarded as elements of $\BB_{2m+1}(\omega)$ over $\CC[\om]$.
As we point out below (see Remark~\ref{rem:allcoef}\tss(i)),
all coefficients are easily expressible in terms of the constant terms
of the polynomials in $\tau$.
Therefore, we will only be concerned with the constant terms modified by the scalar
defined in \eqref{gam}:
\beql{phimde}
\phi_m\coloneqq\ga_m(\om) \q{m} s^{(m)} \f_1 \dots \f_m \q{m}\ts 1,
\eeq
assuming $\tau\ts 1=0$.
We will keep the notation introduced in Sec.~\ref{subsec:mr}, and
for a partition of $\la\vdash m$ of length $\ell=\ell(\la)$,
set
\ben
f[-\la]=f[-\la_1]_1 \dots f[-\la_{\ell}]_{\ell};
\een
cf. \eqref{fla}.

\begin{proposition}\label{prop:integ}
The element \eqref{phimde} is given by the formula
\beql{phimint}
\phi^{}_{m}=\sum_{\la\ts\vdash m,\  \ell(\la)\ts \text{even}}\ts \Yc_{m,\ell}(\om-1)\ts
c^{}_{\la}\ts \q{\ell} h^{(\ell)} f[-\la]\tss \q{m}.
\eeq
\end{proposition}

\bpf
Expand each $\f_a$ as $\tau + f[-1]_a$ and use the second relation in \eqref{taufco}
to move all copies of $\tau$ to the right.
It follows that $\phi_m$ is a linear combination of terms of the form
\beql{proex}
    \gamma_m(\omega) \q{m} s^{(m)} f[-r_1]_{a_1} \dots f[-r_\ell]_{a_\ell} \q{m},
\eeq
where $1 \leq a_1 < \dots < a_\ell \leq m$ and the $r_i$ are positive integers with
$r_1 + \dots + r_\ell = m$.
Now, by \eqref{frafsbabs} for $r,s>0$, we have
\begin{align}
f[-r]_af[-s]_b = f[-s]_b f[-r]_a + \vp_{ab}f[-r-s]_b - f[-r-s]_b \vp_{ab}.
\non
\end{align}
Applying this to each term, we can swap adjacent $f[-r_i]_{a_i}$ until the $r_i$ are weakly decreasing. Since $s^{(m)}\vp_{ab} = \vp_{ab}s^{(m)} = s^{(m)}$, the additional terms with $\vp_{ab}$ cancel out by cyclic property~1. Then, using $s^{(m)}s_{ab} = s_{ab}s^{(m)} = s^{(m)}$ together with
$s_{ab}f[-r]_a = f[-r]_bs_{ab}$ and cyclic property~1, we can insert appropriate permutations $s_{ab}$ after $s^{(m)}$ to change the indices of the $f[-r]_a$ factors to $1$ up to $\ell$, in order.
It follows that
\ben
\phi_m=\sum_{\la\ts\vdash m}\ts
c^{}_{\la}\ts \ga_m(\om)\q{m} s^{(m)} f[-\la]\tss \q{m}.
\een
for some nonnegative integers $c_\lambda$.
By the same argument as in \cite[Theorem~2.1]{m:ss},
we find that $c_\lambda$ is the number of permutations in $\mathfrak{S}_m$ of cycle type $\lambda$.

Furthermore, by repeated application of Lemma~\ref{lem:lemma132}, we have
\ben
\phi_m=\sum_{\la\ts\vdash m}\ts\Yc_{m,\ell}(\om-1)\ts
c^{}_{\la}\ts \ga_{\ell}(\om)\q{\ell} s^{(\ell)} f[-\la]\tss \q{m}.
\een

Now we need a lemma
which
relates these expressions to similar expressions involving $h^{(\ell)}$.

\begin{lemma}\label{lem:lkmequal}
Let $\lambda \vdash m$ be a partition of length $\ell$. Then
\begin{align}
\gamma_\ell(\omega) \q{\ell} s^{(\ell)} f[-\lambda] \q{m}
=
\q{\ell} h^{(\ell)} f[-\lambda] \q{m}.
\non
\end{align}
\end{lemma}

\begin{proof}
From \eqref{eq:gammash}, using the embedding $\B_\ell(\omega) \hookrightarrow \BB_{2m+1}(\omega)$,  we have
\begin{align}
h^{(\ell)} - \gamma_\ell(\omega) s^{(\ell)} \in J_m^{(\ell)} \subseteq \BB_{2m+1},
\non
\end{align}
where $J_m^{(\ell)}$ is spanned by elements $d+d^{t_a}$, where $d$ is a Brauer diagram in $\B_\ell(\omega) \subseteq \BB_{2m+1}(\omega)$ and $1 \leq a \leq \ell$. It thus suffices to show that
\beql{qellze}
\q{\ell} \left(d+d^{t_a}\right) f[-\lambda] \q{\ell} = 0
\eeq
for all $d$.
 We now consider four cases, depending on the diagrams $d$ and $d^{t_a}$. Diagrams with a shaded area are used to represent Brauer diagrams where the nodes within the shaded area are connected by strings that lie within the shaded area; within each case, each diagram with a shaded area has the same arrangement of strings within the shaded area.

\vspace{3mm}

\noindent \textbf{Case (a):} In this case,
\begin{align}
d= d^{t_a} = \;
\begin{tikzpicture}[baseline={([yshift=4.5mm]current bounding box.south)},scale=0.55]
{
\fill[gray!30] (0,2)--(1.75,2) arc(180:360:0.25 and 0.27)--(6,2)--(6,0)--(2.25,0) arc(0:180:0.25 and 0.27)--(0,0)--cycle;
\draw (2,0)--(2,2);
\draw[very thick] (-0.5,0)--(6.5,0);
\draw[very thick] (-0.5,2)--(6.5,2);
\draw (2,2) node[font=\scriptsize,anchor=south]{$a$};
}
\end{tikzpicture}\;  \in \B_\ell(\omega).
\non
\end{align}
\noindent
Noting that $d$ commutes with $s_{ab}$ and $\epsilon_{ab}$ for any $b>\ell$, we have
\begin{multline}
\epsilon_{a\, m+a} \left(d + d^{t_a} \right) f[-\lambda] \epsilon_{a\, m+a}
= 2 \epsilon_{a\, m+a} d\ts f[-\lambda_1]_1 \dots f[-\lambda_\ell]_\ell \epsilon_{a\, m+a}
\non\\[0.4em]
{}= 2 d\ts f[-\lambda_1]_1 \dots f[-\lambda_{a-1}]_{a-1}
\epsilon_{a\, m+a} f[-\lambda_a]_a \epsilon_{a\, m+a} f[-\lambda_{a+1}]_{a+1} \dots f[-\lambda_\ell]_\ell
\end{multline}
which is zero due to the relation $\epsilon_{xy} f[r]_x \epsilon_{xy} = 0$.

\vspace{3mm}
\noindent \textbf{Case (b):} In this case,
\begin{align}
d+d^{t_a} &=\; \begin{tikzpicture}[baseline={([yshift=4.5mm]current bounding box.south)},scale=0.55]
{
\fill[gray!30] (0,2)--(1.75,2) arc(180:360:0.25 and 0.27) --(4.75,2) arc(180:360:0.25 and 0.27)
--(6,2)--(6,0) --(4.25,0) arc(0:180:0.25 and 0.27)--(2.25,0) arc(0:180:0.25 and 0.27)--(0,0)--cycle;
\draw (2,0) to[out=90,in=270] (5,2);
\draw (4,0) to[out=90,in=270] (2,2);
\draw[very thick] (-0.5,0)--(6.5,0);
\draw[very thick] (-0.5,2)--(6.5,2);
\draw (2,2) node[font=\scriptsize,anchor=south]{$a$};
\draw (4,2) node[font=\scriptsize,anchor=south]{$b$};
\draw (5,2) node[font=\scriptsize,anchor=south]{$c$};
}
\end{tikzpicture}
\; +\;
\begin{tikzpicture}[baseline={([yshift=4.5mm]current bounding box.south)},scale=0.55]
{
\fill[gray!30] (0,2)--(1.75,2) arc(180:360:0.25 and 0.27) --(4.75,2) arc(180:360:0.25 and 0.27)
--(6,2)--(6,0) --(4.25,0) arc(0:180:0.25 and 0.27)--(2.25,0) arc(0:180:0.25 and 0.27)--(0,0)--cycle;
\draw (2,0) to[out=90,in=90] (4,0);
\draw (2,2) to[out=270,in=270] (5,2);
\draw[very thick] (-0.5,0)--(6.5,0);
\draw[very thick] (-0.5,2)--(6.5,2);
\draw (2,2) node[font=\scriptsize,anchor=south]{$a$};
\draw (4,2) node[font=\scriptsize,anchor=south]{$b$};
\draw (5,2) node[font=\scriptsize,anchor=south]{$c$};
}
\end{tikzpicture}
\non\\[3mm]
&= \; \begin{tikzpicture}[baseline={([yshift=9.5mm]current bounding box.south)},scale=0.55]
{
\fill[gray!30] (0,2)--(1.75,2) arc(180:360:0.25 and 0.27) --(4.75,2) arc(180:360:0.25 and 0.27)
--(6,2)--(6,0) --(4.25,0) arc(0:180:0.25 and 0.27)--(2.25,0) arc(0:180:0.25 and 0.27)--(0,0)--cycle;
\draw (2,0) -- (2,2);
\draw (4,0) to[out=90,in=270] (5,2);
\begin{scope}[shift={(0,-2)}]
\draw (0,0)--(0,2);
\draw (1.75,0)--(1.75,2);
\draw (2.25,0)--(2.25,2);
\draw (3.75,0)--(3.75,2);
\draw (4.25,0)--(4.25,2);
\draw (6,0)--(6,2);
\draw (2,0) to[out=90,in=270] (4,2);
\draw (4,0) to[out=90,in=270] (2,2);
\foreach \y in {0.2,0,-0.2}{
\node at (0.875+\y,1) {.};
}
\foreach \y in {0.2,0,-0.2}{
\node at (3+\y,0.5) {.};
}
\foreach \y in {0.2,0,-0.2}{
\node at (5.125+\y,1) {.};
}
\end{scope}
\draw[very thick] (-0.5,0)--(6.5,0);
\draw[very thick] (-0.5,2)--(6.5,2);
\draw[very thick] (-0.5,-2)--(6.5,-2);
\draw (2,2) node[font=\scriptsize,anchor=south]{$a$};
\draw (4,2) node[font=\scriptsize,anchor=south]{$b$};
\draw (5,2) node[font=\scriptsize,anchor=south]{$c$};
}
\end{tikzpicture}
\; +\;
\begin{tikzpicture}[baseline={([yshift=9.5mm]current bounding box.south)},scale=0.55]
{
\fill[gray!30] (0,2)--(1.75,2) arc(180:360:0.25 and 0.27) --(4.75,2) arc(180:360:0.25 and 0.27)
--(6,2)--(6,0) --(4.25,0) arc(0:180:0.25 and 0.27)--(2.25,0) arc(0:180:0.25 and 0.27)--(0,0)--cycle;
\draw (2,0)--(2,2);
\draw (4,0) to[out=90,in=270] (5,2);
\begin{scope}[shift={(0,-2)}]
\draw (0,0)--(0,2);
\draw (1.75,0)--(1.75,2);
\draw (2.25,0)--(2.25,2);
\draw (3.75,0)--(3.75,2);
\draw (4.25,0)--(4.25,2);
\draw (6,0)--(6,2);
\draw (2,0) to[out=90,in=90] (4,0);
\draw (2,2) to[out=270,in=270] (4,2);
\foreach \y in {0.2,0,-0.2}{
\node at (0.875+\y,1) {.};
}
\foreach \y in {0.2,0,-0.2}{
\node at (3+\y,1) {.};
}
\foreach \y in {0.2,0,-0.2}{
\node at (5.125+\y,1) {.};
}
\end{scope}
\draw[very thick] (-0.5,0)--(6.5,0);
\draw[very thick] (-0.5,2)--(6.5,2);
\draw[very thick] (-0.5,-2)--(6.5,-2);
\draw (2,2) node[font=\scriptsize,anchor=south]{$a$};
\draw (4,2) node[font=\scriptsize,anchor=south]{$b$};
\draw (5,2) node[font=\scriptsize,anchor=south]{$c$};
}
\end{tikzpicture}
= \left(s_{ab} + \epsilon_{ab} \right) \tilde{d},
\non
\end{align}
where
\begin{align}
\tilde{d}
&= \; \begin{tikzpicture}[baseline={([yshift=4.5mm]current bounding box.south)},scale=0.55]
{
\fill[gray!30] (0,2)--(1.75,2) arc(180:360:0.25 and 0.27) --(4.75,2) arc(180:360:0.25 and 0.27)
--(6,2)--(6,0) --(4.25,0) arc(0:180:0.25 and 0.27)--(2.25,0) arc(0:180:0.25 and 0.27)--(0,0)--cycle;
\draw (2,0)--(2,2);
\draw (4,0) to[out=90,in=270] (5,2);
\draw[very thick] (-0.5,0)--(6.5,0);
\draw[very thick] (-0.5,2)--(6.5,2);
\draw (2,2) node[font=\scriptsize,anchor=south]{$a$};
\draw (4,2) node[font=\scriptsize,anchor=south]{$b$};
\draw (5,2) node[font=\scriptsize,anchor=south]{$c$};
}
\end{tikzpicture}\; .
\non
\end{align}
Note that we do not specify which diagram is which, so this case includes the cases where $d$ is either of these two diagrams.
Observe that $\tilde{d}$ commutes with $s_{a\, m+a}$. From $s_{xy} \epsilon_{xy} = \epsilon_{xy}s_{xy} = \epsilon_{xy}$ and $\ep_{xy} s_{ay} = \ep_{xy}\ep_{ax}$ we have
\begin{align}
\q{\ell} = \q{\ell} s_{b\, m+b} = s_{b\, m+b} \q{\ell}
\Fand \q{\ell}(s_{ab} + \epsilon_{ab})=\q{\ell}(s_{a\ts m+b} + \epsilon_{a\ts m+b}).
\label{qelleq}
\end{align}
Also, from $f[r]_x \epsilon_{xy} = - f[r]_y \epsilon_{xy}$, it follows that
\begin{align}
f[-\lambda]\ts \q{\ell}
&= - f[-\lambda_1]_1 \dots f[-\lambda_b]_{m+b} \dots f[-\lambda_\ell]_\ell\ts \q{\ell}.
\non
\end{align}
Using these properties, we find
\begin{multline}
\q{\ell} \left(d + d^{t_a} \right) f[-\lambda] \q{\ell}
= \q{\ell} \left(s_{ab} + \epsilon_{ab} \right) \tilde{d}\, f[-\lambda] \q{\ell}
\non\\[0.4em]
{}= -\q{\ell} \left(s_{a\, m+b} + \epsilon_{a\, m+b} \right) \tilde{d}\,
f[-\lambda_1]_1 \dots f[-\lambda_b]_{m+b} \dots f[-\lambda_\ell]_\ell \q{\ell}.
\end{multline}
Using the first relation in \eqref{qelleq}, we find that this coincides with
$-\q{\ell} \left(d + d^{t_a} \right) f[-\lambda] \q{\ell}$
and hence is equal to zero. Note that this proof does not
rely on the ordering of $a$, $b$ and $c$, and holds even if $b=c$.

\vspace{3mm}
\noindent \textbf{Case (c):} In this case,
\begin{align}
d+d^{t_a} &=\; \begin{tikzpicture}[baseline={([yshift=4.5mm]current bounding box.south)},scale=0.55]
{
\fill[gray!30] (0,2)--(1.25,2) arc(180:360:0.25 and 0.27) --(6,2)--(6,0) -- (4.75,0) arc(0:180:0.25 and 0.27)--(3.25,0) arc(0:180:0.25 and 0.27) --(1.75,0) arc(0:180:0.25 and 0.27)--(0,0)--cycle;
\draw (3,0) to[out=90,in=270] (1.5,2);
\draw (1.5,0) to[out=90,in=90] (4.5,0);
\draw[very thick] (-0.5,0)--(6.5,0);
\draw[very thick] (-0.5,2)--(6.5,2);
\draw (1.5,2) node[font=\scriptsize,anchor=south]{$a$};
\draw (3,2) node[font=\scriptsize,anchor=south]{$b$};
\draw (4.5,2) node[font=\scriptsize,anchor=south]{$c$};
}
\end{tikzpicture}
\; +\;
\begin{tikzpicture}[baseline={([yshift=4.5mm]current bounding box.south)},scale=0.55]
{
\fill[gray!30] (0,2)--(1.25,2) arc(180:360:0.25 and 0.27) --(6,2)--(6,0) -- (4.75,0) arc(0:180:0.25 and 0.27)--(3.25,0) arc(0:180:0.25 and 0.27) --(1.75,0) arc(0:180:0.25 and 0.27)--(0,0)--cycle;
\draw (4.5,0) to[out=90,in=270] (1.5,2);
\draw (1.5,0) to[out=90,in=90] (3,0);
\draw[very thick] (-0.5,0)--(6.5,0);
\draw[very thick] (-0.5,2)--(6.5,2);
\draw (1.5,2) node[font=\scriptsize,anchor=south]{$a$};
\draw (3,2) node[font=\scriptsize,anchor=south]{$b$};
\draw (4.5,2) node[font=\scriptsize,anchor=south]{$c$};
}
\end{tikzpicture}
= \left(s_{ab} + \epsilon_{ab} \right) \tilde{d},
\non
\end{align}
where
\begin{align}
\tilde{d}
&= \; \begin{tikzpicture}[baseline={([yshift=4.5mm]current bounding box.south)},scale=0.55]
{
\fill[gray!30] (0,2)--(1.25,2) arc(180:360:0.25 and 0.27) --(6,2)--(6,0) -- (4.75,0) arc(0:180:0.25 and 0.27)--(3.25,0) arc(0:180:0.25 and 0.27) --(1.75,0) arc(0:180:0.25 and 0.27)--(0,0)--cycle;
\draw (3,0) to[out=90,in=90] (4.5,0);
\draw (1.5,0) -- (1.5,2);
\draw[very thick] (-0.5,0)--(6.5,0);
\draw[very thick] (-0.5,2)--(6.5,2);
\draw (1.5,2) node[font=\scriptsize,anchor=south]{$a$};
\draw (3,2) node[font=\scriptsize,anchor=south]{$b$};
\draw (4.5,2) node[font=\scriptsize,anchor=south]{$c$};
}
\end{tikzpicture}\; .
\non
\end{align}

Then, since this $\tilde{d}$ also commutes with $s_{a\, m+a}$, the computation from the previous case holds here, so $\q{\ell} \left(d + d^{t_a} \right) f[-\lambda] \q{\ell} = 0$.

\vspace{3mm}
\noindent \textbf{Case (d):} In this case, both $d$ and $\tilde{d}$ are obtained
from the respective diagrams in Case (c) by reflection about a horizontal line, so that $d+d^{t_a}=  \tilde{d} (s_{ab} + \epsilon_{ab} )$.
Again, $s_{a\, m+a}$ commutes with $\tilde{d}$. The computation to show $\q{\ell} \left(d+ d^{t_a}\right) f[-\lambda] \q{\ell} = 0$ in this case is similar to the one in the previous two cases: simply replace $\left(s_{ab} + \epsilon_{ab} \right) \tilde{d}$ with $\tilde{d} \left(s_{ab} + \epsilon_{ab}\right)$, and $\left(s_{a\, m+b} + \epsilon_{a\, m+b} \right) \tilde{d}$ with $\tilde{d} \left(s_{a\, m+b}
+ \epsilon_{a\, m+b} \right)$, throughout.
\end{proof}

Note that if $\ell$ is odd, then $h^{(\ell)}$ belongs to the subspace $J_m^{(\ell)}$.
This follows by writing
\ben
2\tss h^{(\ell)}=\frac{1}{\ell!}\ts\sum_{s\in\Sym_\ell} (s+s^{-1})
\een
and using the telescoping sum
\ben
s+s^{-1}=(s+s^{t_1})-(s^{t_1}+s^{t_1t_2})+(s^{t_1t_2}+s^{t_1t_2t_3})
-\dots+(s^{t_1\dots t_{\ell-1}}+s^{t_1\dots t_{\ell}}).
\een
Therefore,
$\q{\ell} h^{(\ell)} f[-\lambda] \q{m}=0$ by \eqref{qellze}
(in particular, $\phi_1=0$). This observation and
Lemma~\ref{lem:lkmequal}
complete the proof of the proposition.
\epf

We can now use Propositions~\ref{prop:thm816first}, \ref{prop:thm816second} and \ref{prop:integ}
to complete the proof of Theorem~\ref{thm:osp}. The element $\phi^{}_{m}$
given in \eqref{phimint} belongs to the algebra $\BB_{2m+1}(\omega)$
defined over $\CC[\om]$. Consider this algebra at the critical level by taking
its quotient $\BB_{2m+1}(\omega)_{\cri}$ over the ideal generated
by the central element $\om+K-2$.
We conclude that
both $f[0]_0\ts \phi_m$ and $f[1]_0\ts \phi_m$ belong to
the left ideal of the algebra $\BB_{2m+1}(\omega)_{\cri}$ over $\CC(\om)$ generated by
the subspaces $f[r]_a\tss \B_{2m+1}(\om)$ for $r\geq 0$
and $a=0,\dots,m$.
To be able to evaluate $\om$ at $M-2n$ we need to verify that the elements of the left ideal
are defined over  $\CC[\om]$. This is clear for $f[0]_0\ts \phi_m$
from Propositions~\ref{prop:thm816first} and \ref{prop:integ},
whereas
the formula of Proposition~\ref{prop:thm816second} implies
the relation in the algebra $\BB_{2m+1}(\omega)_{\cri}$:
\ben
\bal
f[1]_0\ts \phi_m
{}={}& - m\ts\ga_m(\om)\ts \q{m} \vp_{0m}s^{(m)}  \f_1 \dots \f_{m-1} f[0]_m \q{m}\ts 1\\[0.5em]
&+ m\ts\ga_m(\om)\ts\q{m} s^{(m)}  \f_1 \dots \f_{m-1} \q{m} f[0]_0 \ts 1\\[0.5em]
& + m\ts\ga_m(\om)\ts \q{m} s^{(m)} \vp_{0m} \f_1 \dots \f_{m-1} f[0]_m \q{m}\ts 1\\[0.5em]
& + \ga_m(\om)\ts\q{m}\ts s^{(m)} \f_1 \dots \f_m \q{m} f[1]_0\ts 1.
 \eal
\een
By cyclic property 2 and relations
\ben
\vp_{0m}\ts f[0]_m=f[0]_0\ts s_{0m}+\ep_{0m}\ts f[0]_0\Fand
f[0]_m\ts \vp_{0m}=s_{0m}\ts f[0]_0+ f[0]_0\ts\ep_{0m},
\een
the sum of the first and third terms
on the right hand side equals
\ben
m\ts \ga_m(\om)\ts\q{m}s^{(m)}  \f_1 \dots \f_{m-1} f[0]_0 \vp_{0m}\q{m}\ts 1
- m\ts \ga_m(\om)\ts\q{m} s^{(m)} \f_1 \dots \f_{m-1}\vp_{0m}\ts f[0]_0\ts \q{m}\ts 1.
\een
With this replacement, the argument used in the proof of Proposition~\ref{prop:integ}
applies to all four terms in the expansion of $f[1]_0\ts \phi_m$
with obvious modifications
to show that all of them are defined over $\CC[\om]$.
For example, expanding the second of the new terms, we can write it as a $\CC$-linear combination
of the expressions analogous to \eqref{proex}:
\ben
\gamma_m(\omega) \q{m} s^{(m)} f[-r_1]_{a_1} \dots f[-r_{\ell-1}]_{a_{\ell-1}} \vp_{0m}\ts f[0]_0\ts\q{m},
\een
where $1 \leq a_1 < \dots < a_{\ell-1} \leq m-1$ and the $r_i$ are positive integers with
$r_1 + \dots + r_{\ell-1} = m-1$. We can insert appropriate permutations $s_{ab}$
after $s^{(m)}$ to change the indices to write
such an expression as
\ben
\gamma_m(\omega) \q{m} s^{(m)} f[-r_1]_1 \dots f[-r_{\ell-1}]_{\ell-1} \vp_{0\tss\ell}\ts f[0]_0\ts\q{m}.
\een
Now apply Lemma~\ref{lem:lemma132} repeatedly to find that it equals
\ben
\gamma_\ell(\omega) \q{\ell} s^{(\ell)} f[-r_1]_1 \dots f[-r_{\ell-1}]_{\ell-1} \vp_{0\tss\ell}\ts f[0]_0\ts\q{m}
\een
times an element of $\CC[\om]$. Finally, Lemma~\ref{lem:lkmequal} applies in the same
form (the role of $f[-\la_{\ell}]_{\ell}$ is played by $\vp_{0\tss\ell}$)
allowing us to replace $\gamma_\ell(\omega)s^{(\ell)}$ with $h^{(\ell)}$
to produce an element of the algebra $\BB_{2m+1}(\omega)_{\cri}$ over $\CC[\om]$.

By applying the homomorphism $\rho$ defined in \eqref{homrho}, we find
that the image $\rho(\phi_m)$ coincides with $\Phi_m \ts \Q{m}$, where
$\Phi_m$ is defined in \eqref{phimliesup}.
Under the evaluation $\om=M-2n$,
the central element $\om+K-2$ becomes $K+h^{\vee}$ so that by
the annihilation properties of $\phi_m$ in $\BB_{2m+1}(\omega)_{\cri}$,
$\Phi_m$
belongs to the Feigin--Frenkel centre $\z(\wh\osp_{M|2n})$,
thus completing the proof of Theorem~\ref{thm:osp}.

\medskip

\begin{remark}\label{rem:allcoef}
\quad (i)
It was pointed out in \cite[Remark~2.5]{m:ss} that if $M=0$ or $n=0$, then
all coefficients of the polynomial in $\tau$ appearing in
\eqref{defphi} coincide with the Segal--Sugawara vectors $\Phi_k$ for certain $k$, up
to a constant factor. A similar property is shared by the polynomial
$\ga_m(\om) \q{m} s^{(m)} \f_1 \dots \f_m \q{m}$ or, more generally, by
\beql{genssk}
\ga_k(\om) \q{k} s^{(k)} \f_1 \dots \f_k\ts \q{m}=\psi^{}_{k0}\ts \tau^k+\dots+\psi^{}_{kk},
\eeq
where we keep $m$ fixed, while $1 \leq k \leq m$. Namely, we have the relations
\ben
    \psi^{}_{ka}= \binom{\omega + k - 2}{k-a}\ts \psi^{}_{aa}.
\een
They are easily verified by noting that for any $u\in\CC$ the map $\tau\mapsto u+\tau$
extends to an automorphism of the algebra $\BB_{2m+1}(\omega)$. We apply it to
both sides of \eqref{genssk} and compare the coefficients of $u^{k-a}\tau^0$.
\par
\quad (ii) We can apply the evaluation homomorphism $\ev_z:\osp_{M|2n}[t^{-1}]\to \osp_{M|2n}$
which takes $t$ to a nonzero complex number $z$, to the Segal--Sugawara vectors $\Phi_m$.
As with the non-super case considered in \cite[Sec.~6.5]{m:so}, the images $\ev_z(\Phi_m)$
belong to the centre of the universal enveloping algebra $\U(\osp_{M|2n})$;
cf.~\cite[Thm~3.16]{lw:sw}.
\end{remark}

\subsection*{Acknowledgements}
Our work is supported by the Australian Research Council, grant DP240101572.

\bigskip
\bigskip

\small

\noindent
School of Mathematics and Statistics\newline
University of Sydney,
NSW 2006, Australia\\[0.5em]
{\tt alexander.molev@sydney.edu.au}\\
{\tt madeline.nurcombe@sydney.edu.au}

\end{document}